\documentclass[preprint,12pt,numbers,sort&compress]{elsarticle}
\usepackage{soul}
\usepackage{psfrag}
\usepackage{a4wide}
\usepackage{rotating}
\usepackage{color}
\usepackage{amssymb}
\usepackage{amsmath}
\usepackage{amsthm}
\usepackage{verbatim}
\usepackage{cancel}

\newcommand{\f}[1]{\mathbf{#1}}
\newcommand{\ab}[1]{\boldsymbol{#1}}
\def\bfm#1{\boldsymbol{#1}}
\newcommand{\bb}[1]{\bfm{#1}}

\newcommand{\R}{\mathbb R}

\newcommand{\V}{\mathcal{V}}

\newcommand{\W}{\mathcal{W}}
\newcommand{\glob}{y}
\newcommand{\loc}{\zeta}
\newcommand{\LL}{i_0}
\newcommand{\RR}{i_1}
\newcommand{\g}{f}
\newcommand{\gC}{f}
\newcommand{\sm}{4}
\newcommand{\Side}{\tau}
\newcommand{\dd}{\partial}
\newcommand{\ff}{f^{\bfm{\Xi}^{(i)}}}

\DeclareMathOperator{\Span}{span}

\theoremstyle{definition}
\newtheorem{ex}{Example}

\newproof{pf}{proof}
\advance\textheight by 0.4cm
\advance\topmargin by -0.2cm

\definecolor{gold}{rgb}{1,0.7,0}
\definecolor{dred}{rgb}{0.92,0,0}
\definecolor{dgreen}{rgb}{0,0.6,0}

\def\new{\color{black}}

\begin{document}

\begin{frontmatter}

\title{Isogeometric collocation for solving the biharmonic equation over planar multi-patch domains}

\cortext[cor]{Corresponding author}

\author[vil]{Mario Kapl}
\ead{m.kapl@fh-kaernten.at}

\author[slo1]{Alja\v z Kosma\v c\corref{cor}}
\ead{aljaz.kosmac@iam.upr.si}

\author[slo1]{Vito Vitrih}
\ead{vito.vitrih@upr.si}

\address[vil]{ADMiRE Research Center, Carinthia University of Applied Sciences, Villach, Austria}

\address[slo1]{IAM and FAMNIT, University of Primorska, Koper, Slovenia}

\begin{abstract}
We present an isogeometric collocation method for solving the biharmonic equation over planar bilinearly parameterized multi-patch domains. The developed approach
is based on the use of the globally $C^4$-smooth isogeometric spline space~\cite{KaVi20b} to approximate the solution of the considered partial differential equation, and proposes as collocation points two different choices, namely on the one hand the Greville points and on the other hand the so-called superconvergent points. Several examples demonstrate the potential of our collocation method for solving the biharmonic equation over planar multi-patch domains, and numerically study the convergence behavior of the two types of collocation points with respect to the $L^2$-norm as well as to equivalents of the $H^s$-seminorms for $1 \leq s \leq 4$.
{\new 
In the studied case of spline degree $p=9$, the numerical results indicate in case of the Greville points a convergence of order~$\mathcal{O}(h^{p-3})$ independent of the considered (semi)norm, and show in case of the superconvergent points an improved convergence of order~$\mathcal{O}(h^{p-2})$ for all (semi)norms except for the equivalent of the $H^4$-seminorm, where the order $\mathcal{O}(h^{p-3})$ is anyway optimal.
}
\end{abstract}


\end{frontmatter}

\section{Introduction}

Isogeometric Analysis (IgA) is a numerical framework for solving a partial differential equation (PDE) by employing the same spline or NURBS functions to represent the geometry of the computational domain as well as the numerical solution of the considered PDE~\cite{CottrellBook,HuCoBa04}.
The most common technique in IgA is to transform the original problem (also called strong form of the PDE) into its weak form and to solve the weak form by means of the Galerkin method. The benefit of this approach is that the numerical solution of the weak form can be of (much) lower regularity than the solution of the given strong form. However, the method is based on the evaluation of integrals of (rational) functions of high degree,
{\new which requires high computational costs of the matrix assembly routines. Furthermore, consistency, robustness and order of convergence are important issues that have to be considered carefully for each numerical integration technique. An important progress has been achieved in the last years (e.g. \cite{CaSaTa2017,GiJu2022,HuReSa2010,JuMaPeRu2016,MaJu2015,PaJuGi2020}), but numerical integration still remains an issue.}
Another approach in IgA, which we will follow in this work, is to solve the given strong form of the PDE directly via the collocation method. In this way no integration is needed {\new and matrix assembly is much faster
}, but the method requires spaces of higher smoothness like $C^2$ and $C^4$-smooth isogeometric spline functions for second and fourth order PDEs, respectively.

Considering planar multi-patch domains with possible extraordinary vertices like in this paper, globally $C^s$-smooth ($s \geq 1$) isogeometric spline functions are usually constructed by means of the concept of geometric continuity of multi-patch surfaces~\cite{Pe02} using the fact that an isogeometric function is $C^s$-smooth if and only if the associated multi-patch graph surface is $G^s$-smooth~\cite{Pe15,KaViJu15}. While for~$s=1$ several techniques exist, see the survey articles~\cite{HuSaTaTo21,KaSaTa19b} for more details, there are just a few methods available for~$s \geq 2$, namely the constructions~\cite{KaVi17a,KaVi17b,KaVi17c,KaVi19a,ToSpHiHu16} for $s=2$ and the approach \cite{KaVi20b} for an arbitrary $s \geq 1$, {\new which includes the case~$s=4$ needed for solving the biharmonic equation in strong form via collocation.} 

{\new
Besides multi-patch quadrangular domains, triangulations can be used to generate $C^s$-smooth spline spaces over complex domains. The book \cite{LaSch07} gives an overview of different techniques to model such smooth spline spaces, and provides a detailed bibliography on this topic. Some recent constructions of $C^4$-smooth spline spaces on triangulations using polynomial macro-element spaces and Powell-Sabin splits
can be find e.g.~in \cite{Speleers2013, 
Groselj2016}.
}

The problem of isogeometric collocation has been firstly explored in~\cite{IsoCollocMethods2010} for solving 
second order PDEs on a single patch. Thereby, an important issue in isogeometric collocation is the selection of the collocation points, since it affects the convergence behavior of the numerical solution. E.g. for second order problems in case of odd spline degrees~$p$, the application of Greville points (abscissae) leads to a convergence of order $\mathcal{O}(h^{p-1})$ with respect to $h$-refinement independent of considering the $L^2$, $H^1$ or $H^2$-error~\cite{IsoCollocMethods2010}, but the use of so-called superconvergent (Cauchy Galerkin) points~\cite{SuperConvergent2015,GomezLorenzisVariationalCollocation,MonSanTam2017} improves this convergence to an order of $\mathcal{O}(h^{p+1})$ ($\mathcal{O}(h^{p})$ for~\cite{GomezLorenzisVariationalCollocation}) and $O(h^{p})$ for the $L^2$ and $H^1$-error, respectively. Also for solving fourth order PDEs, Greville and superconvergent points are mainly used. So far, the numerical comparison of the convergence behavior between the two different types of collocation points have been performed just for the $L^2$~norm~\cite{GomezLorenzisVariationalCollocation}, where for odd spline degrees~$p$ the superconvergent points show again a better convergence behavior, namely a convergence of order~$\mathcal{O}(h^{p-2})$ compared to~$\mathcal{O}(h^{p-3})$ for the Greville points. 

While in case of one-patch domains, isogeometric collocation has been already extensively used for solving second and fourth order problems such as the Poisson's equation, e.g.~\cite{SuperConvergent2015,IsoCollocMethods2010,FahrenderLorenzisGomez2018,GomezLorenzisVariationalCollocation,MonSanTam2017,CostComparison2013}, problems of linear and nonlinear elasticity, elastostatistics and elastodynamics, e.g.~\cite{SuperConvergent2015,Elatostatics2012,EvansHughesReali2018,FahrenderLorenzisGomez2018,GomezLorenzisVariationalCollocation,FaMoReHuLo20,JiaAnitescuZhangRabczuk2019,Reali2015,MoFaLoEvHuRe21,CostComparison2013}, acoustics problems, e.g.~\cite{AyViAnAt20,AtCaAnKh22, ZaPa21}, the Reissner-Mindlin shell
problem, e.g.~\cite{EnzoKiendlLorenzis2017}, beam problems, e.g.~\cite{RealiGomez2015,ShearTimoshenko2012,EnzoKiendlLorenzis2019}, the biharmonic equation/Kirchhoff plate problem, e.g.~\cite{TsplinesIgC2016,RealiGomez2015, GomezLorenzisVariationalCollocation, Maurin2018}, the Kirchhoff-Love shell problem, e.g.~\cite{Maurin2018} and phase-field models, e.g.~\cite{GomezRealiSangali2014,Schillinger2015}, in case of multi-patch domains, the number of existing work is small, and is further limited to the solving of second order problems, namely to the Poisson's equation~\cite{KaVi20} as well as to elasticity problems~\cite{Elatostatics2012,JiaAnitescuZhangRabczuk2019}. In the case of multi-patch domains, two different strategies have been used for solving the considered second order PDEs. The technique~\cite{KaVi20} employs the globally $C^2$-smooth isogeometric spline space~\cite{KaVi19a} as discretization space for the Poisson's equation, which allows the direct solving of the strong form of the second order PDE. In contrast, the methods~\cite{Elatostatics2012,JiaAnitescuZhangRabczuk2019} employ spline spaces, which are not globally $C^2$-smooth having a lower regularity at the patch interfaces, and require therefore a special treatment at the patch interface in the collocation procedure. 

This paper extends the multi-patch isogeometric collocation method~\cite{KaVi20} for second order problems to the solving of fourth order PDEs demonstrated on the basis of the biharmonic equation. The presented approach uses the globally $C^4$-smooth isogeometric spline space~\cite{KaVi20b} to solve the biharmonic equation over bilinearly parameterized multi-patch domains in a strong form, i.e. to compute a $C^4$-smooth numerical solution of the PDE. Two different choices of collocation points, given by the Greville points as well as by a particular family of superconvergent points, are numerically studied for the case of the odd spline degree~$p=9$ and inner patch regularity~$r=4$, which is the configuration with the lowest possible degree for the $C^4$-smooth spline space~\cite{KaVi20b}. The convergence behavior with respect to $h$-refinement is not only investigated for the $L^2$-norm as in~\cite{GomezLorenzisVariationalCollocation,Maurin2018}, but also for (equivalents to) the $H^s$-seminorms for $1 \leq s \leq 4$. The numerical results indicate in case of the Greville points a convergence of order~$\mathcal{O}(h^{p-3})$ independent of the considered norm, and show in case of the superconvergent points an improved convergence of order~$\mathcal{O}(h^{p-2})$ for all (semi)norms except for the equivalent of the $H^4$-seminorm, where the order $\mathcal{O}(h^{p-3})$ is anyway optimal. Furthermore, it is demonstrated that the rates of convergence for a multi-patch domain coincide with the ones for a one-patch domain.  

The remainder of the paper is organized as follows. Section~\ref{sec:problem_statement} introduces the general framework for solving the biharmonic equation over planar bilinear multi-patch domains via isogeometric collocation. For this purpose, we describe the used multi-patch structure of the domain, explain the general concept of $C^4$-smooth isogeometric spline functions over planar bilinear multi-patch domains, and present the multi-patch isogeometric collocation approach for solving the biharmonic equation. In Section~\ref{sec:collocation}, we determine the final collocation method by specifying on the one hand the used $C^4$-smooth isogeometric discretization space, namely the $C^4$-smooth spline space~\cite{KaVi20b}, to represent the $C^4$-smooth approximation of the solution of the biharmonic equation, and by introducing on the other hand two possible choices of collocation points. Section~\ref{section_Numerical_examples} presents several numerical examples, which study the convergence behavior under $h$-refinement with respect to the $L^2$-norm as well as to (equivalents of) the $H^s$-seminorms for $1 \leq s \leq 4$, and which demonstrate the power of our collocation method to solve the biharmonic equation over bilinearly parameterized multi-patch domains. Finally, we conclude the paper in Section~\ref{sec:Conclusion}.

\section{Multi-patch isogeometric collocation for the biharmonic equation} \label{sec:problem_statement}

We will describe the general framework for performing isogeometric collocation for solving the biharmonic equation over planar bilinearly parameterized multi-patch domains. This includes the presentation of the considered multi-patch configuration of the domain, the introduction of the concept of $C^4$-smooth isogeometric spline functions over these multi-patch domains as well as the development of the multi-patch isogeometric collocation technique for solving the biharmonic equation. 

\subsection{The multi-patch configuration of the domain} \label{subsec:multipatch}

\begin{figure}[bth]
    \centering
    \includegraphics[scale=0.27]{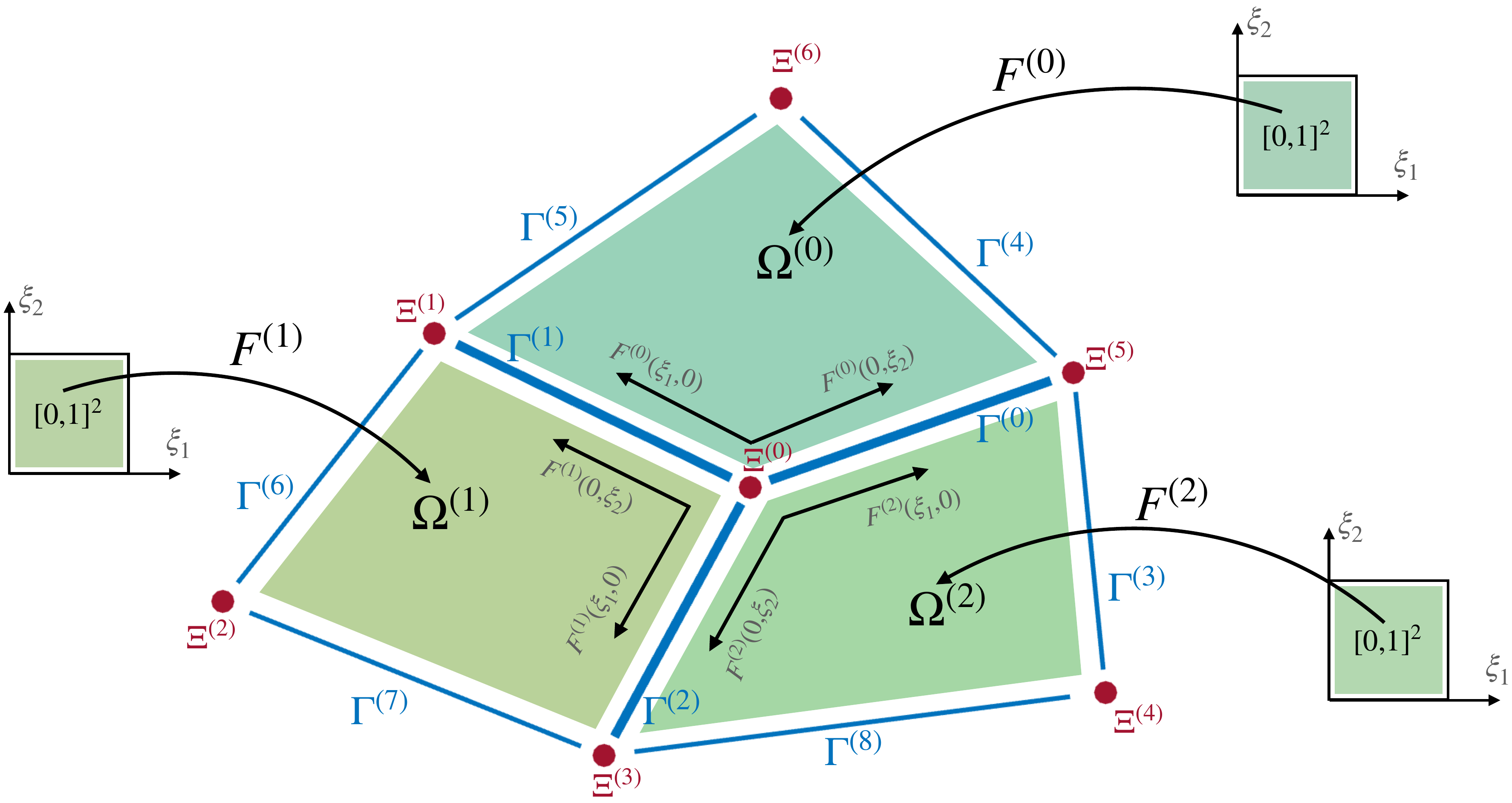}
    \caption{An example of a multi-patch domain $\overline{\Omega}$  with the patches $\Omega^{(i)}$ and their associated geometry mappings $\ab{F}^{(i)}$ (green), with the edges $\Gamma^{(i)}$ (blue) and with the vertices $\bfm{\Xi}^{(i)}$ (red).}
    \label{fig:multipatchCase}
\end{figure}

Let $\Omega \subset \R^2$ be an open and connected planar domain, whose closure $\overline{\Omega}$ is the multi-patch domain~$\overline{\Omega} = \cup_{i \in \mathcal{I}_{\Omega}} \overline{\Omega^{(i)}}$ with open quadrilateral patches~$\Omega^{(i)}$, $i \in \mathcal{I}_{\Omega}$, which can be further represented as the disjoint union of the open patches~$\Omega^{(i)}$, $i \in \mathcal{I}_{\Omega}$, of open edges~$\Gamma^{(i)}$, $i \in \mathcal{I}_{\Gamma}$, and of vertices~$\bfm{\Xi}^{(i)}$, $i \in \mathcal{I}_{\Xi}$, i.e.
\begin{equation*} 
\displaystyle
\overline{\Omega} = \bigcup_{i \in \mathcal{I}_{\Omega}} \Omega^{(i)}  \; \dot{\cup}  \bigcup_{i \in \mathcal{I}_{\Gamma}} \Gamma^{(i)} \; \dot{\cup} \bigcup_{i \in \mathcal{I}_{\Xi}} \bfm{\Xi}^{(i)},
\end{equation*}
where $\dot{\cup}$ denotes the disjoint union of sets, and where $\mathcal{I}_{\Omega}$, $\mathcal{I}_{\Gamma}$ and $\mathcal{I}_{\Xi}$ are the index sets of the indices of the patches~$\Omega^{(i)}$, edges~$\Gamma^{(i)}$ and vertices $\bfm{\Xi}^{(i)}$, respectively. The index sets $\mathcal{I}_{\Gamma}$ and $\mathcal{I}_{\Xi}$ are further divided into $\mathcal{I}_{\Gamma} = \mathcal{I}_{\Gamma}^I \dot{\cup} \mathcal{I}_{\Gamma}^B$ and $\mathcal{I}_{\Xi} = \mathcal{I}_{\Xi}^I \dot{\cup} \mathcal{I}_{\Xi}^B$, where $\mathcal{I}_{\Gamma}^I$ and $\mathcal{I}_{\Gamma}^B$ collect the indices of the inner and boundary edges~$\Gamma^{(i)}$, respectively, and where $\mathcal{I}_{\Xi}^I$ and $\mathcal{I}_{\Xi}^B$ contain the indices of the inner and boundary vertices~$\bfm{\Xi}^{(i)}$, respectively. We will assume that the closures of any two patches~$\Omega^{(i_0)}$ and $\Omega^{(i_1)}$, $i_0,i_1 \in \mathcal{I}_{\Omega}$, have either an empty intersection, i.e. $\overline{\Omega^{(i_0)}} \cap \overline{\Omega^{(i_1)}} = \emptyset$, possess exactly one common vertex, i.e. $\overline{\Omega^{(i_0)}} \cap \overline{\Omega^{(i_1)}} = \bfm{\Xi}^{(i)}$ for an $i \in \mathcal{I}_{\Xi}$, or share the 
closure of a common inner edge, i.e. $\overline{\Omega^{(i_0)}} \cap \overline{\Omega^{(i_1)}} = \overline{\Gamma^{(i)}}$ for an $i \in \mathcal{I}_{\Gamma}^I$. Furthermore, we will assume that each patch $\overline{\Omega^{(i)}}$ is parameterized by a bilinear, bijective and regular geometry mapping~$\ab{F}^{(i)}$,
\begin{align*}
 \ab{F}^{(i)}: [0,1]^{2}  \rightarrow \R^{2}, \quad 
 \bb{\xi} =(\xi_1,\xi_2) \mapsto
 \ab{F}^{(i)}(\bb{\xi}) = \ab{F}^{(i)}(\xi_1,\xi_2), \quad i \in \mathcal{I}_{\Omega},
\end{align*}
such that $\overline{\Omega^{(i)}} = \ab{F}^{(i)}([0,1]^{2})$, cf. Fig.~\ref{fig:multipatchCase}.

\subsection{The concept of $C^4$-smooth isogeometric functions over bilinear multi-patch domains}

\begin{figure}
    \centering
    \includegraphics[scale=0.23]{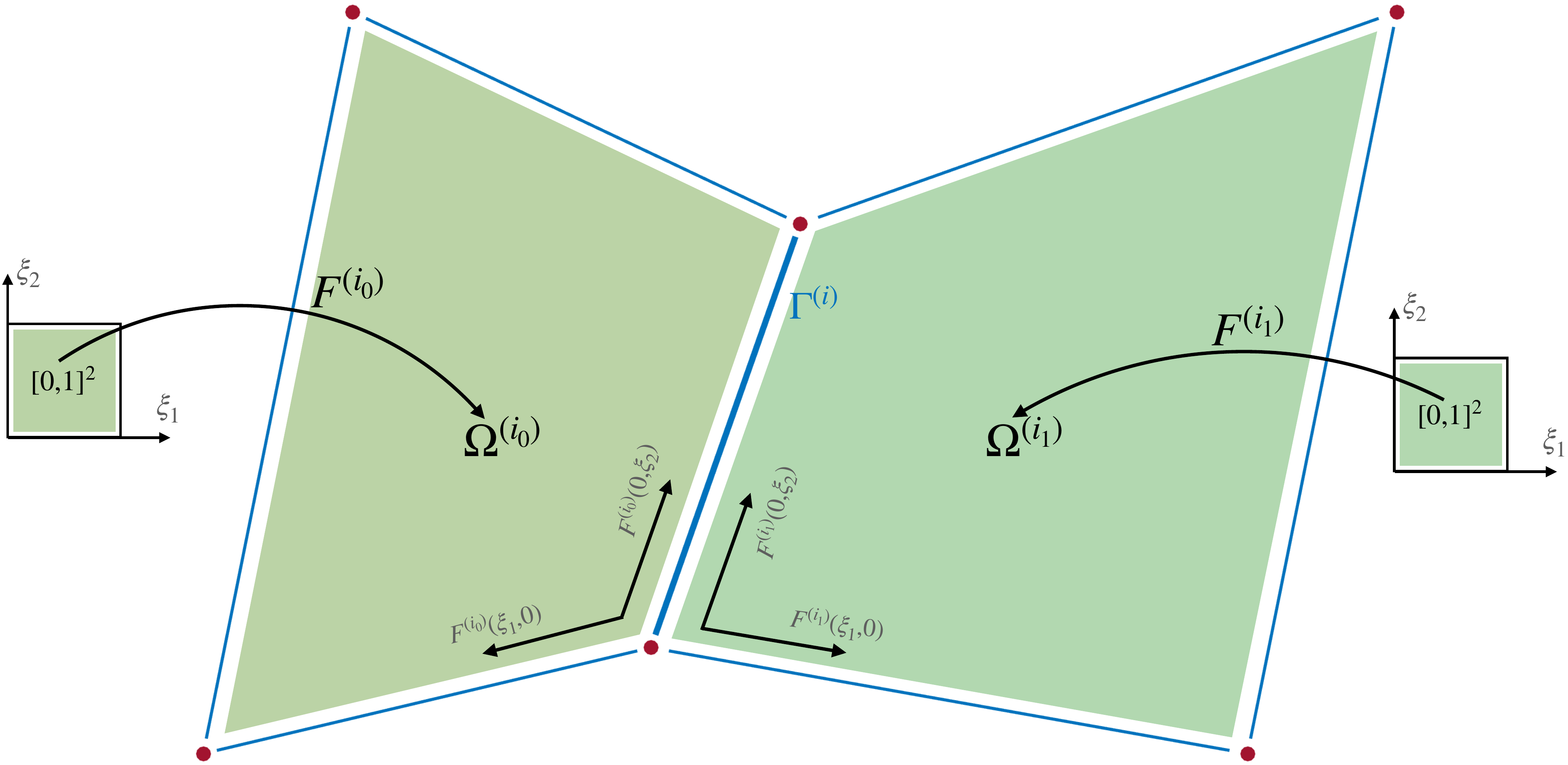}
    \caption{The parameterization of the two-patch domain $\overline{\Omega^{(i_0)} \cup \Omega^{(i_1)}}$ with the common interface $\Gamma^{(i)}$ and corresponding geometry mappings $\ab{F}^{(i_0)}$ and $\ab{F}^{(i_1)}$.}
    \label{fig:twopatchCase}
\end{figure}

Let $\mathcal{S}_h^{p,r}([0,1])$ be the univariate spline space of degree~$p$, regularity~$r$ and mesh size~$h=\frac{1}{k+1}$ defined on the unit interval~$[0,1]$ and constructed from the uniform open knot vector
\begin{equation*}  
(t_0^{p,r},\ldots,t_{2p+k(p-r)+1}^{p,r})=
(\underbrace{0,\ldots,0}_{(p+1)-\mbox{\scriptsize times}},
\underbrace{\textstyle \frac{1}{k+1},\ldots,\frac{1}{k+1}}_{(p-r) - \mbox{\scriptsize times}},\ldots, 
\underbrace{\textstyle \frac{k}{k+1},\ldots ,\frac{k}{k+1}}_{(p-r) - \mbox{\scriptsize times}},
\underbrace{1,\ldots,1}_{(p+1)-\mbox{\scriptsize times}}),
\end{equation*}
where $k$ is the number of different inner knots. Furthermore, let $\mathcal{S}_h^{\ab{p},\ab{r}}([0,1]^2)$ be the tensor-product spline space $\mathcal{S}_h^{p,r}([0,1]) \otimes 
\mathcal{S}_h^{p,r}([0,1])$ on the unit-square~$[0,1]^2$, and let $N_{j}^{p,r}$ and $N_{j_1,j_2}^{\ab{p},\ab{r}}=N_{j_1}^{p,r}N_{j_2}^{p,r}$ with $j,j_1,j_2=0,1,\ldots,n-1$, and $n= \dim \mathcal{S}_h^{p,r}([0,1]) = p+1+k(p-r)$, the B-splines of the spline spaces~$\mathcal{S}_h^{p,r}([0,1])$ and $\mathcal{S}_h^{\ab{p},\ab{r}}([0,1]^2)$, respectively. 
Since the geometry mappings~$\ab{F}^{(i)}$, $i \in \mathcal{I}_{\Omega}$, are bilinearly parameterized, we trivially have that
\[
\ab{F}^{(i)} \in \mathcal{S}_{h}^{\ab{p},\ab{r}}([0,1]^2) \times \mathcal{S}^{\ab{p},\ab{r}}_{h}([0,1]^2).
\]
The space of $C^4$-smooth isogeometric spline functions over the multi-patch domain~$\overline{\Omega}$ 
is given as
\begin{equation*}
\V^4 = \V \cap C^4(\overline{\Omega}),
\end{equation*}
where
\begin{equation*}
\V = \left\{ \phi \in L^2(\overline{\Omega}): \; \phi |_{\overline{\Omega}^{(i)}} \in {\mathcal{S}_{h}^{\ab{p},\ab{r}}([0,1]^{2})} \circ (\ab{F}^{(i)})^{-1} , \; 
i \in \mathcal{I}_{\Omega}   \right\}
\end{equation*}
is the space of isogeometric spline functions over the multi-patch domain~$\overline{\Omega}$. An isogeometric spline function $\phi \in \V$ belongs to the space~$\mathcal{V}^4$ if and only if for any two neighboring patches~$\Omega^{(i_0)}$ and $\Omega^{(i_1)}$, $i_0,i_1 \in \mathcal{I}_{\Omega}$, with the common inner edge~$\overline{\Gamma^{(i)}} = \overline{\Omega^{(i_0)}} \cap 
\overline{\Omega^{(i_1)}}$, $i \in \mathcal{I}_{\Gamma}^I$, the associated graph surface patches
 $ [ \ab{F}^{(i_0)}, f^{(i_0)} ]^T$
and
$[ \ab{F}^{(i_1)}, f^{(i_1)}]^T$, 
with the spline functions $f^{(\Side)}=\phi \circ \ab{F}^{(\Side)}$, $\Side \in \{\LL,\RR \}$, are $G^4$-smooth across their common interface, cf.~\cite{Pe15,KaViJu15}. 
An equivalent condition to the $G^4$-smoothness of the two neighboring graph surface patches is that 
the two associated spline functions~$f^{(i_0)}$ and $f^{(i_1)}$ satisfy for all $\xi \in [0,1]$
\begin{equation}   \label{eq:gC}
 \g_\ell^{(i,\LL)}(\xi) = \g_\ell^{(i,\RR)}(\xi) =: \gC^{(i)}_\ell(\xi), \quad \ell =0,1,\ldots,4,
 \end{equation}
 with
 \begin{equation*}   \label{eq:gC2}
 \g_\ell^{(i,\Side)}(\xi) = \left(\alpha^{(i,\Side)}(\xi)\right)^{-\ell}\, \partial_1^\ell \g^{(\Side)}(0,\xi) - \sum_{i=j}^{\ell-1} {\ell \choose j} 
 \left(\frac{\beta^{(i,\Side)}(\xi)}{\alpha^{(i,\Side)}(\xi)}\right)^{\ell-j}  \dd^{\ell-j} \gC_j^{(i)}(\xi) ,\quad \Side\in \{\LL,\RR\},
 \end{equation*}
where the functions $\alpha^{(i,\Side)},\beta^{(i,\Side)}, \Side \in \{\LL,\RR \}$, {\new are called gluing functions, and} are linear polynomials defined as 
\begin{equation*}  \label{eq:alphaLRbar}
 \alpha^{(i,\Side)}(\xi) = \lambda \det J \ab{F}^{(\Side)}(0,\xi), 
 \quad \beta^{(i,\Side)}(\xi) = \frac{\partial_1 \ab{F}^{(\Side)}(0,\xi) \cdot \partial_2\ab{F}^{(\Side)}(0,\xi)}
 {||\partial_2 \ab{F}^{(\Side)}(0,\xi)||^{2}}, \quad \lambda > 0,
\end{equation*}
and $J \ab{F}^{(\Side)}$, $\Side \in \{\LL,\RR \}$, are the Jacobians of $\ab{F}^{(\Side)}$, assuming without loss of generality that the two corresponding geometry mappings~$\ab{F}^{(i_0)}$ and $\ab{F}^{(i_1)}$ are parameterized as shown in Fig.~\ref{fig:twopatchCase}, cf.~\cite[Theorem 2]{KaVi20b}. Moreover, we have $\alpha^{(i,i_0)} < 0$ and $\alpha^{(i,i_1)} > 0$ due to the regularity of the geometry mappings~$\ab{F}^{(i_0)}$ and $\ab{F}^{(i_1)}$, and can select $\lambda$ e.g. by minimizing 
$
 || \alpha^{(i,\LL)}+1 ||^2_{L^2} + || \alpha^{(i,\RR)}-1 ||^2_{L^2},
$
cf. \cite{KaVi19a}. Using the $C^4$-smoothness conditions~\eqref{eq:gC}, the space~$\V^4$ can be now described as
\begin{equation*}
\V^4 = \left\{ \begin{array}{ll} \phi \in L^2(\overline{\Omega}): \; & \phi |_{\overline{\Omega}^{(i)}} \in {\mathcal{S}_{h}^{\ab{p},\ab{r}}([0,1]^{2})} \circ (\ab{F}^{(i)})^{-1} , \; 
i \in \mathcal{I}_{\Omega}, \mbox{ and } \\ & \g_\ell^{(i,\LL)}(\xi) = \g_\ell^{(i,\RR)}(\xi), \; \xi \in [0,1], \, \ell=0,1,\ldots, 4,\; i \in \mathcal{I}_{\Gamma}^I \end{array}   \right\}.
\end{equation*}

\subsection{Isogeometric collocation for the biharmonic equation} \label{subsec:biharmonic}

The goal is to find  $u:\overline{\Omega} \to \R$, $u \in C^4(\overline{\Omega})$, which solves the biharmonic equation
\begin{align} \label{eq:biharmonic}
\triangle^2 u (\bfm{x}) & =   g(\bfm{x}), \quad \; \bfm{x} \in \overline{\Omega},  \nonumber \\
u (\bfm{x}) & = g_1(\bfm{x}) ,   \quad \bfm{x} \in \partial \Omega,  \\
\partial_{\bfm{n}} u (\bfm{x}) & = g_2(\bfm{x}) ,   \quad \bfm{x} \in \partial \Omega, \nonumber
\end{align}
{\new where $g:\overline{\Omega} \to \R$ and $g_1, g_2: \partial \Omega \to \R$ are sufficiently smooth functions.}
We will use the concept of isogeometric collocation to compute a $C^4$-smooth approximation $u_h \in \mathcal{W}^4$ of the solution~$u$, where $\mathcal{W}^4 \subseteq \V^4$ is a suitable discretization space. 
For this we need a set of global collocation points~$\bfm{\glob}_j$, $j \in \mathcal{J}$, which are further separated into inner collocation points $\bfm{\glob}_j^I$, $j \in \mathcal{J}_I$, and boundary collocation points $\bfm{\glob}_j^B$,  $j \in \mathcal{J}_B$. Inserting these points into \eqref{eq:biharmonic}, we obtain 
\begin{align} \label{eq:biharmonicPoints}
\triangle^2 u (\bfm{\glob}_j^I) & =   g(\bfm{\glob}_j^I), \quad \;\, j \in \mathcal{J}_I,\nonumber \\
u (\bfm{\glob}_j^B) & = g_1(\bfm{\glob}_j^B) ,   \quad j \in \mathcal{J}_B,  \\
\partial_{\bfm{n}} u (\bfm{\glob}_j^B) & = g_2(\bfm{\glob}_j^B) ,   \quad j \in \mathcal{J}_B.  \nonumber  \nonumber 
\end{align}
To use the isogeometric approach for solving problem~$\eqref{eq:biharmonicPoints}$, we have first to express each global collocation point~$\bfm{\glob}_{j}$ with respect to local coordinates via
$$
\bfm{\loc}^{I}_{j} = \left( \bfm{F}^{(\iota_j)}\right)^{-1} \left( \bfm{\glob}_j^I \right) , \;  j \in \mathcal{J}_I, \quad \mbox{ and } \quad
\bfm{\loc}^{B}_{j} = \left( \bfm{F}^{(\iota_j)}\right)^{-1} \left( \bfm{\glob}_j^B \right) , \;  j \in \mathcal{J}_B,
$$
where 
\begin{equation} \label{eq:patch_index_selection}  
\iota_j := \min \{ i \in \mathcal{I}_\Omega, \; \bfm{\glob}_j \in \overline{\Omega^{(i)}} \}.
\end{equation}
Once we have the local collocation points $\bfm{\loc}^{I}_{j}$ and $\bfm{\loc}^{B}_{j}$, equations \eqref{eq:biharmonicPoints} can be transformed (\cite{BaDe15}) to
\begin{align}  \label{eq:collocationSystemLocal}
\frac{1}{\left| \det J \ab{F}^{(\iota_j)} \left( \bfm{\loc}^{I}_{j}  \right) \right|}  \nabla \circ  \hspace{-0.05cm} \Bigg( N^{(\iota_j)}(\bb{\xi}) 
\nabla \Bigg(  \frac{1}{\left| \det J \ab{F}^{(\iota_j)}  \left( \bfm{\xi} \right) \right|}  & \nabla \circ \left( N^{(\iota_j)}(\bb{\xi}) 
\nabla \left(u \left( \bfm{F}^{(\iota_j)}(\bb{\xi}) \right) \right) \right) \hspace{-0.13cm} \bigg) \hspace{-0.05cm} \bigg) {\Bigg|}_{\bfm{\xi} =\bfm{\loc}^{I}_{j} } \nonumber\\
&
=
g\left(\ab{F}^{(\iota_j)}\left(\bfm{\loc}^{I}_{j}\right)\right),\;\, \mbox{ }j \in \mathcal{J}_{I}, \nonumber \\
&  \\[-0.4cm]
\normalsize
u \left( \bfm{F}^{(\iota_j)}  \left(\bfm{\loc}^{B}_{j} \right)\right)  & = g_1\left(\ab{F}^{(\iota_j)}\left(\bfm{\loc}^{B}_{j}\right)\right),
 \mbox{ }j \in \mathcal{J}_{B}, \nonumber \\
 \left\langle\bfm{n},(J\bfm{F}^{(\iota_j)}\left(\bfm{\loc}^{B}_{j}\right))^{-1}\nabla u\left( \bfm{F}^{(\iota_j)}  \left(\bfm{\loc}^{B}_{j} \right)\right) \right\rangle  & = g_2\left(\ab{F}^{(\iota_j)}\left(\bfm{\loc}^{B}_{j}\right)\right),
 \mbox{ }j \in \mathcal{J}_{B}, \nonumber
\end{align}
where 
\begin{equation*}
N^{(\iota_j)}(\bb{\xi}) = \left|\det J \ab{F}^{(\iota_j)}(\bb{\xi})\right|  \left(J \ab{F}^{(\iota_j)}(\bb{\xi})\right)^{-T}  \left(J \ab{F}^{(\iota_j)}(\bb{\xi})\right)^{-1} .
\end{equation*}
This leads to a linear system for the unknown coefficients $c_i$ of the approximation $u_h = \sum_{i \in \mathcal{I}} c_i \phi_i$, where $\mathcal{I} = \{0,1, \ldots, \dim \mathcal{W}^4-1 \}$, and $\{\phi_i\}_{i\in \mathcal{I}}$ is the basis of $\mathcal{W}^4$. In the following section, we will finally determine our isogeometric multi-patch collocation method for solving the biharmonic equation by specifying the used $C^4$-smooth isogeometric discretization space~$\W^4$ and by presenting two possible sets of collocation points.

\section{The specific setting in the isogeometric collocation} \label{sec:collocation}

We will present the specific setting for the $C^4$-smooth isogeometric discretization space $\W^4$ and for the set of collocation points in the isogeometric multi-patch collocation approach introduced in the previous section.  

\subsection{The $C^4$-smooth isogeometric discretization space~$\W^4$} \label{C4_isogeometric_splinespace}

Instead of the entire $C^4$-smooth isogeometric spline space~$\V^4$, which is difficult to construct, we will use as discretization space in the isogeometric collocation method the simpler $C^4$-smooth subspace
\[
\W^4 = \{ \phi \in \V^4 : \; \gC^{(i)}_\ell \in \mathcal{S}_h^{p-\ell,r+4-\ell}([0,1]) , \, \ell=0,1,\ldots, 4,  \, i \in \mathcal{I}_{\Gamma}^{I}\} \subseteq \, \V^4,
\]
which has been introduced in~\cite{KaVi20b} for an arbitrary smoothness~$s \geq 1$. The construction of the space~$\W^4$, which will be summarized\footnote{For full details of the construction of the $C^4$-smooth spline space~$\W^4$, we refer to~\cite{KaVi20b}.} below, will require a spline degree~$p \geq 9$, an inner patch regularity~$r \geq 4$ with $4 \leq r \leq p-5$ as well as a mesh size~$h \leq\frac{p-r-4}{12-r+1}$. In this paper, we will restrict ourselves to the case of the lowest possible spline degree~$p$ and regularity~$r$, namely to $p=9$ and $r=4$, which further implies $h \leq \frac{1}{9}$. 
 
The $C^4$-smooth isogeometric spline space~$\W^4$ is generated as the direct sum of smaller subspaces corresponding to the single patches~$\Omega^{(i)}$, $i \in \mathcal{I}_\Omega$, edges~$\Gamma^{(i)}$, $i \in \mathcal{I}_\Gamma$, and vertices $\bfm{\Xi}^{(i)}$, $i \in \mathcal{I}_\Xi$, i.e., 
\begin{equation*} \label{eq:W}
    \mathcal{W}^4 =  \left(\bigoplus_{i \in \mathcal{I}_\Omega} \mathcal{W}^4_{\Omega^{(i)}}\right) \oplus 
    \left(\bigoplus_{i \in \mathcal{I}_\Gamma} \mathcal{W}^4_{\Gamma^{(i)}}\right) 
    \oplus  
    \left(\bigoplus_{i \in \mathcal{I}_\Xi} \mathcal{W}^4_{\bfm{\Xi}^{(i)}}\right).
\end{equation*}
The single subspaces $\W^4_{\Omega^{(i)}}$, $\W^4_{\Gamma^{(i)}}$ and  $\W^4_{\bfm{\Xi}^{(i)}}$ are constructed as the span of corresponding basis functions as demonstrated in the following paragraphs. 

\paragraph{The patch subspace~$\W^4_{\Omega^{(i)}}$} 

For each patch~$\Omega^{(i)}$, $i \in \mathcal{I}_{\Omega}$, the patch subspace $\mathcal{W}^4_{\Omega^{(i)}}$ is defined as
\begin{equation*} \label{eq:spaceW0hOmega}
\mathcal{W}^4_{\Omega^{(i)}} =  \Span \left\{ \phi_{\Omega^{(i)};j_1,j_2} |\;  j_1,j_2=5,6,\ldots,n-6  \right\},
\end{equation*}
with the isogeometric functions
\begin{equation}  \label{eq:PhiOmega2}
{\phi}_{\Omega^{(i)};j_1,j_2}(\bfm{x})  = 
\begin{cases}
   (N_{j_1,j_2}^{\bfm{9},\bfm{4}}\circ (\ab{F}^{(i)})^{-1})(\bfm{x}) \;
\mbox{ if }\f \, \bfm{x} \in \overline{\Omega^{(i)}},
\\ 0 \quad \mbox{ if }\f \, \bfm{x} \in \overline{\Omega} \backslash \overline{\Omega^{(i)}}.
\end{cases} 
\end{equation} 

\paragraph{The edge subspace~$\W^4_{\Gamma^{(i)}}$} 

It has to be distinguished between the case of a boundary and of an inner edge~$\Gamma^{(i)}$. In case of a boundary edge $\Gamma^{(i)}$, $i \in \mathcal{I}_{\Gamma}^{B}$, with $\Gamma^{(i)} \subseteq \overline{\Omega^{(i_0)}}$, 
$i_0 \in \mathcal{I}_{\Omega}$, assuming without loss of generality that the boundary edge $\Gamma^{(i)}$ is given by
$\ab{F}^{(i_0)}(\{0 \} \times (0,1)) $, the edge subspace $\mathcal{W}^4_{\Gamma^{(i)}}$ is defined as
\begin{align*} 
\mathcal{W}^4_{\Gamma^{(i)}} =&  \Span \left\{ \phi_{\Omega^{(i_0)}; j_1,j_2} |\; j_2=9-j_1,\ldots,n+j_1-10, \; j_1=0,1,\ldots,4 \right\} ,
\end{align*}
where the isogeometric functions~$\phi_{\Omega^{(i_0)}; j_1,j_2}$ are given as in~\eqref{eq:PhiOmega2}. 

In case of an inner edge $\Gamma^{(i)}$, $i \in \mathcal{I}_{\Gamma}^I$ with $\Gamma^{(i)} \subseteq  \overline{\Omega^{(i_0)}} \cap \overline{\Omega^{(i_1)}}$, $i_0,i_1 \in \mathcal{I}_{\Omega}$, assuming without loss of generality that the two associated geometry mappings~$\ab{F}^{(i_0)}$ and $\ab{F}^{(i_1)}$ are parameterized as shown in Fig.~\ref{fig:twopatchCase}, the edge subspace $\mathcal{W}^4_{\Gamma^{(i)}}$ is given as
\begin{equation*} \label{eq:spaceW0hGamma}
\mathcal{W}^4_{\Gamma^{(i)}} =  \Span \left\{\phi_{\Gamma^{(i)};j_1,j_2}| \;\; j_2=9-j_1,
\ldots,k(j_1+1),\; j_1=0,1,\ldots,4 \right\},
\end{equation*} 
where the isogeometric functions $\phi_{\Gamma^{(i)};j_1,j_2} $ 
have the form
\begin{equation}  \label{eq:defphiGamma}
    \phi_{\Gamma^{(i)};j_1,j_2} (\bfm{x}) = 
\begin{cases}
  (\g_{\Gamma^{(i)}; j_1,j_2}^{(i_0)} \circ (\ab{F}^{(i_0)})^{-1})(\bfm{x}) \;
\mbox{ if }\f \, \bfm{x} \in \overline{\Omega^{(i_0)}},
\\[0.15cm] 
  (\g_{\Gamma^{(i)}; j_1,j_2}^{(i_1)} \circ (\ab{F}^{(i_1)})^{-1} )(\bfm{x}) \;
\mbox{ if }\f \, \bfm{x} \in \overline{\Omega^{(i_1)}}, 
\\[0.15cm] 
  0\quad {\rm otherwise},
\end{cases}
\end{equation}
with
\begin{align}
 \g_{\Gamma^{(i)}; 0,j_2}^{(\Side)}(\xi_1,\xi_2) & = \left(\beta^{(\Side)}(\xi_2)\right)^4 \partial^4 N_{j_2}^{9,8}(\xi_2) M_4 (\xi_1) +\left(\beta^{(\Side)}(\xi_2)\right)^3 \partial^3 N_{j_2}^{9,8}(\xi_2) M_3 (\xi_1) + \nonumber\\ 
 & \quad \left(\beta^{(\Side)}(\xi_2)\right)^2 \partial^2 N_{j_2}^{9,8}(\xi_2) M_2 (\xi_1) +  
 \beta^{(\Side)}(\xi_2) \partial N_{j_2}^{9,8}(\xi_2) M_1 (\xi_1) + N_{j_2}^{9,8}(\xi_2) M_0 (\xi_1), \nonumber\\
 \g_{\Gamma^{(i)}; 1,j_2}^{(\Side)}(\xi_1,\xi_2) & =  \frac{9}{h} \left( 4 \left(\beta^{(\Side)}(\xi_2)\right)^3 \alpha^{(\Side)}(\xi_2) \partial^3 N_{j_2}^{8,7}(\xi_2) M_4 (\xi_1) + 3 \left(\beta^{(\Side)}(\xi_2)\right)^2 \alpha^{(\Side)}(\xi_2) \partial^2 N_{j_2}^{8,7}(\xi_2)\cdot \right. \nonumber\\ 
 & \quad \left. M_3 (\xi_1) + 2 \beta^{(\Side)}(\xi_2)\alpha^{(\Side)}(\xi_2) \partial N_{j_2}^{8,7}(\xi_2) M_2 (\xi_1) +  
 \alpha^{(\Side)}(\xi_2) N_{j_2}^{8,7}(\xi_2) M_1(\xi_1)
 \right), \nonumber
\end{align} 
\vskip-2.5em
\begin{align}  \label{eq:basisFunctionsGenericG}
 \g_{\Gamma^{(i)}; 2,j_2}^{(\Side)}(\xi_1,\xi_2) & = \frac{72}{h^2}\left( 6 \left(\beta^{(\Side)}(\xi_2)\right)^2 \left(\alpha^{(\Side)}(\xi_2)\right)^{2} \partial^2 N_{j_2}^{7,6}(\xi_2) M_4 (\xi_1) + \right. \\ 
 &  \quad \left. 3 \beta^{(\Side)}(\xi_2) \left(\alpha^{(\Side)}(\xi_2)\right)^{2} \partial N_{j_2}^{7,6}(\xi_2) M_3 (\xi_1) + \left(\alpha^{(\Side)}(\xi_2)\right)^{2} N_{j_2}^{7,6}(\xi_2) M_2 (\xi_1) \right), \nonumber\\
 \g_{\Gamma^{(i)}; 3,j_2}^{(\Side)}(\xi_1,\xi_2) & = \frac{504}{h^3} \left( 4 \beta^{(\Side)}(\xi_2) \left(\alpha^{(\Side)}(\xi_2)\right)^{3} \partial N_{j_2}^{6,5}(\xi_2) M_4 (\xi_1) + \left(\alpha^{(\Side)}(\xi_2)\right)^{3} N_{j_2}^{6,5}(\xi_2) M_3 (\xi_1) \right),\nonumber\\
 \g_{\Gamma^{(i)}; 4,j_2}^{(\Side)}(\xi_1,\xi_2) & = \frac{3024}{h^4} \left(\alpha^{(\Side)}(\xi_2)\right)^{4}  N_{j_2}^{5,4}(\xi_2)  M_4 (\xi_1), \nonumber
\end{align}
and
\begin{equation*}  \label{eq:M012}
  M_i(\xi) = \sum_{j=i}^4 \frac{ {j \choose i} h^i}{\prod_{\ell=0}^{i-1} (9-\ell)} N_j^{9,4}(\xi), \quad i = 0,1,\ldots,4.
\end{equation*}

\paragraph{The vertex subspace~$\W^4_{\bfm{\Xi}^{(i)}}$} 

The construction of the vertex subspace~$\W^4_{\bfm{\Xi}^{(i)}}$ differs whether we have an inner or a boundary vertex~$\bfm{\Xi}^{(i)}$. Let us first consider the case of an inner vertex~$\bfm{\Xi}^{(i)}$, $i \in \mathcal{I}_{\Xi}^I$, with patch valency $v_i$. We assume without loss of generality that all patches $\Omega^{(i_\rho)}$, $\rho=0,1,\ldots,v_i-1$, around the vertex~$\bfm{\Xi}^{(i)}$, i.e. $\bfm{\Xi}^{(i)} = \cap_{\rho=0}^{v_i-1} \overline{\Omega^{(i_\rho)}}$, are parameterized as shown in Fig.~\ref{fig:multipatchCase}, and relabel the common edges $\overline{\Omega^{(i_\rho)}} \cap \overline{\Omega^{(i_{\rho+1})}}$, $\rho=0,1,\ldots,v_i-1$, by $\Gamma^{(i_{\rho+1})}$, where~$\rho$ is taken modulo~$v_i$. The vertex subspace~$\W^4_{\bfm{\Xi}^{(i)}}$ is now defined as
\begin{equation*} 
\mathcal{W}^4_{\bfm{\Xi}^{(i)}} =  \Span \left\{\phi_{\bfm{\Xi}^{(i)};j}| \;\; j=0, \ldots,n^{(i)}-1 \right\},
\end{equation*} 
where the $C^4$-smooth isogeometric functions~$\phi_{\bfm{\Xi}^{(i)};j} $ in the vicinity of the vertex~$\bfm{\Xi}^{(i)}$ are 
given as a linear combination of functions~$\phi_{\Gamma^{(i_{\rho})};j_1,j_2}$, $j_2 = 0,1,\ldots, 
8 -j_1,\; j_1=0,1,\ldots,4$, $\rho = 0,1,\ldots, v_i-1$, coinciding at their common supports in the vicinity of the vertex~$\bfm{\Xi}^{(i)}$, and by subtracting functions~$\phi_{\Omega^{(i_{\rho})}; j_1, j_2}$, $j_1,j_2=0,1,\ldots,4$, $\rho=0,1,\ldots,v_i-1$, which have been added twice. To describe the construction of these functions~$\phi_{\bfm{\Xi}^{(i)};j}$ in more detail, we first consider the isogeometric function $\phi_{\bfm{\Xi}^{(i)}}$ of the form
\begin{equation*}  \label{eq:defphiXi}
  \phi_{\bfm{\Xi}^{(i)}} (\ab{x}) = 
  \begin{cases}
   \left( \ff_{i_\rho} \circ (\ab{F}^{(i_\rho)})^{-1}\right)(\bfm{x}) \;
\mbox{ if }\f \, \bfm{x} \in \overline{\Omega^{(i_\rho)}},\; \rho=0,1,\ldots,{v}_i -1,
\\
0 \quad \mbox{ otherwise,}
\end{cases}
\end{equation*}
where the functions~$ \ff_{i_{\rho}}$ are given as 
\begin{equation*}  \label{eq:g_vertex}
 \ff_{i_{\rho}}(\xi_1,\xi_2)    =   f_{i_\rho}^{\Gamma^{(i_\rho)}}(\xi_1,\xi_2)
  +  f_{i_\rho}^{\Gamma^{(i_{\rho+1})}}(\xi_1,\xi_2) \, - 
 f_{i_\rho}^{\Omega^{(i_\rho)}}(\xi_1,\xi_2), 
\end{equation*}
  with the functions
\begin{align} \label{eq:g_vertex2} 
    f_{i_\rho}^{\Gamma^{(i_{\rho+\tau})}}(\xi_1,\xi_2)  & = \sum_{j_1=0}^\sm \sum_{j_2=0}^{8 -j_1} a^{\Gamma^{(i_{\rho+\tau})}}_{j_1,j_2} \, 
  f_{\Gamma^{(i_{\rho+\tau}); j_1,j_2}}^{(i_\rho)} (\xi_{2-\tau},\xi_{1+\tau}), \quad { a^{\Gamma^{(i_{\rho+\tau})}}_{j_1,j_2} \in \R},
  \quad \tau = 0,1, \nonumber \\[-0.3cm] 
  \\[-0.3cm]
 f_{i_\rho}^{\Omega^{(i_\rho)}}(\xi_1,\xi_2)  & =  \sum_{j_1=0}^{\sm} \sum_{j_2=0}^{\sm } a^{(i_\rho)}_{j_1,j_2} N_{j_1,j_2}^{\ab{p},\ab{r}} 
 (\xi_1,\xi_2), \quad {a_{j_1,j_2}^{(i_\rho)} \in \R }, \nonumber
\end{align}
and with the functions~$f_{\Gamma^{(i_{\rho+\tau}); j_1,j_2}}^{(i_\rho)}$, $\tau=0,1,$ given in~\eqref{eq:basisFunctionsGenericG}.
The isogeometric function~$\phi_{\bfm{\Xi}^{(i)}}$ is $C^4$-smooth on $\overline{\Omega}$ if the coefficients 
$a^{\Gamma^{(i_{\rho+\tau})}}_{j_1,j_2}$, $a_{j_1,j_2}^{(i_\rho)}$ satisfy the equations
\begin{equation} \label{eq:vertex_homogeneous_system}
  \partial_{1}^{\ell_1}  
  \partial_{2}^{\ell_2} 
   \left( f_{i_\rho}^{\Gamma^{(i_{\rho+1})}} - f_{i_\rho}^{\Gamma^{(i_\rho)}} \right)  (\bfm{0}) = 0 \quad \mbox{ and } \quad
    \partial_{1^{}}^{\ell_1}
   \partial_{2^{}}^{\ell_2}
   \left( f_{i_\rho}^{\Gamma^{(i_{\rho+1})}} - f_{i_\rho}^{\Omega^{(i_\rho)}} \right)  (\bfm{0}) = 0, 
\end{equation} 
for $ 0 \leq \ell_1, \ell_2 \leq 4$ and $\rho =0,1, \ldots, v_{i}-1$. The equations~\eqref{eq:vertex_homogeneous_system} form a homogeneous system of linear equations, and a basis of the kernel of this system defines now via the corresponding coefficients $a^{\Gamma^{(i_{\rho+\tau})}}_{j_1,j_2}$, $a_{j_1,j_2}^{(i_\rho)}$ the isogeometric functions $\phi_{\bfm{\Xi}^{(i)};j}$, $j=0,1, \ldots, n^{(i)}-1$, where $n^{(i)}$ is the dimension of the kernel of the homogeneous linear system. In our examples in Section~\ref{section_Numerical_examples}, we will use the concept of minimal determining sets (cf. \cite{KaVi17a,LaSch07}) to determine the basis of the kernel.

Let now $\bfm{\Xi}^{(i)}$, $i \in \mathcal{I}_{\bfm{\Xi}}^{B}$, be a boundary vertex with patch valency~$\nu_i \geq 1$, and let the two boundary edges be labeled as $\Gamma^{(i_0)}$ and $\Gamma^{(i_{v_{i}})}$. In case of $\nu_i \geq 3$, the construction of the space~$\mathcal{W}^{4}_{\bfm{\Xi}^{(i)}}$ works analogously as for an inner vertex, except that for the patches~$\Omega^{(i_0)}$ and $\Omega^{(i_{v_i-1})}$ the functions~$f_{i_0;j_1,j_2}^{\Gamma^{(i_0)}}$ and $f_{i_{v_{i}-1};j_1,j_2}^{\Gamma^{(i_{v_i})}}$ in~\eqref{eq:g_vertex2} are the standard B-splines. In case of $v_i \in \{1,2 \}$, the vertex subspace is directly constructed without solving a homogeneous linear system~\eqref{eq:vertex_homogeneous_system}. For a patch valency $v_i=2$, we assume without loss of generality that the two neighboring patches~$\Omega^{(i_0)}$ and $\Omega^{(i_1)}$, $i_0,i_1 \in \mathcal{I}_{\Omega}$, which contain the vertex~$\bfm{\Xi}^{(i)}$ and possess the common edge~$\overline{\Gamma^{(j_0)}}= \overline{\Omega^{(i_0)}} \cap \overline{\Omega^{(i_1)}}, j_0 \in \mathcal{I}_{\Gamma}$, are parameterized as shown in Fig.~\ref{fig:twopatchCase} and that the vertex~$\bfm{\Xi}^{(i)}$ is further given as 
$\bfm{\Xi}^{(i)} = \ab{F}^{(i_0)}(\ab{0}) = \ab{F}^{(i_1)}(\ab{0})$. The vertex subspace~$\mathcal{W}^{4}_{\bfm{\Xi}^{(i)}}$ is then constructed as
\begin{align*} \label{eq:spaceWXB}
\mathcal{W}^\sm_{\bfm{\Xi}^{(i)}} =&  \Span \left\{ \widetilde{\phi}_{\bfm{\Xi}^{(i)}; j_1,j_2} |\;   
j_1=0,1,\ldots,12, \; j_2 = \begin{cases}
                         0,1,\ldots, 8-j_1 & \mbox{if }j_1 \leq 8 \\
                         0,1,\ldots, 12-j_1 & \mbox{if }j_1 > 8 
                                          \end{cases}
 \; \right \} ,
\end{align*}
with
\begin{equation*} \label{eq:phi_vertex_2}
 \widetilde{\phi}_{\bfm{\Xi}^{(i)};j_1,j_2} (\bfm{x}) = 
    \begin{cases} 
     \phi_{\Gamma^{(j_0)};j_1,j_2}(\ab{x}) & \mbox{if }j_1=0,\ldots ,4\\
     \phi_{\Omega^{(i_{0})};j_1,j_2}(\ab{x}) & \mbox{if }j_1=5,\ldots,8\\
     \phi_{\Omega^{(i_{1})};j_1-4,j_2}(\ab{x}) & \mbox{if }j_1=9,\ldots,12,
    \end{cases}
\end{equation*}
where the functions $\phi_{\Omega^{(i_{\ell})};j_1,j_2}$, $\ell=0,1$, and $\phi_{\Gamma^{(j_0)}_{};j_1,j_2}$ are defined as 
in~\eqref{eq:PhiOmega2} and \eqref{eq:defphiGamma}, respectively. 
For a patch valency $v_i=1$, we assume without loss of generality that $\bfm{\Xi}^{(i)}= \ab{F}^{(i_0)}(\ab{0})$, $i_0 \in \mathcal{I}_{\Omega}$. Then, 
the vertex subspace is constructed as 
$$
\mathcal{W}^4_{\bfm{\Xi}^{(i)}} =  \Span \left\{ \phi_{\Omega^{(i_0)}; j_1,j_2} |\;   
j_1,j_2 =  0,1,\ldots, 8, \; j_1+j_2 \leq 8 \right\} ,
$$
with functions $\phi_{\Omega^{(i_0)}; j_1,j_2}$ given as in~\eqref{eq:PhiOmega2}.

\paragraph{Remark about the condition for the mesh size~$h$} The condition $h \leq\frac{p-r-4}{12-r+1}$, i.e. $h \leq \frac{1}{9}$ for our case $p=9$ and $r=4$, for the mesh size~$h$, ensures the uniform construction of the $C^4$-smooth spline space~$\W^4$ for an arbitrary given planar bilinearly parameterized multi-patch domain~$\overline{\Omega}$. Thereby, the required mesh size~$h$ guarantees that the intersection of two vertex subspaces~$\W^4_{\bfm{\Xi}^{(i)}}$, $i \in \mathcal{I}_{\Xi}$, only contains the zero function. However, this condition for the mesh size~$h$ can be often relaxed. E.g., if the multi-patch domain~$\overline{\Omega}$ possesses no inner vertices and no boundary vertices with a patch valency greater than two, cf. Examples~\ref{ex:onepatchdomain} and \ref{ex:fourpatchdomainHole}, then the $C^4$-smooth space~$\W^4$ technically consists just of patch functions~$\phi_{\Omega^{(i)};j_1,j_2}$ and edge functions~$\phi_{\Gamma^{(i)};j_1,j_2}$. Ensuring that each patch or edge function will be only taken once in case it would belong to more than one subspace, the construction of the $C^4$-smooth space~$\W^4$ will then work for any mesh size~$h \leq 1$. We can proceed similarly if the multi-patch domain~$\overline{\Omega}$ possesses exactly one vertex~$\bfm{\Xi}^{(i)}$, $i \in \mathcal{I}_{\Xi}$, with a patch valency $\nu_i \geq 3$, cf.~Example~\ref{ex:threepatchdomain}. Then we just have additionally to ensure that the patch and edge functions which are involved in the construction of the vertex subspace~$\W^4_{\bfm{\Xi}^{(i)}}$ of the corresponding vertex~$\bfm{\Xi}^{(i)}$ will not be also taken in another subspace. 

\subsection{Selection of collocation points} \label{sec:SelectionCollocPoints}

We will select the global collocation points directly via local collocation points $\bfm{\loc}^{(i)}_{\bfm{j}} = (\loc^{(i)}_{j_1},\loc^{(i)}_{j_2}) \in [0,1]^2 $ for each patch~$\overline{\Omega^{(i)}}$, $i \in \mathcal{I}_{\Omega}$, as $\ab{F}^{(i)}(\bfm{\loc}^{(i)}_{\bfm{j}}) \in \overline{\Omega^{(i)}}$, where $\loc^{(i)}_j \in [0,1]$ are univariate collocation points. For the sake of simplicity, we will use for each patch~$\overline{\Omega^{(i)}}$, $i \in \mathcal{I}_{\Omega}$, just the same local collocation points $\bfm{\loc}_{\bfm{j}} = (\loc_{j_1},\loc_{j_2}) \in [0,1]^2 $ with the univariate collocation points $\loc^{}_j \in [0,1]$. If two or more local collocation points from different patches define the same global collocation point, then we keep just one of these local collocation points, namely the one for the patch~$\overline{\Omega^{(i)}}$ with the lowest patch index~$i$, cf. \eqref{eq:patch_index_selection}, to prevent the possible repetition of global collocation points. 

In the numerical examples in Section~\ref{section_Numerical_examples}, we will study two different choices of univariate collocation points~$\loc^{}_j$, which will then define the local collocation points~$\bfm{\loc}_{\bfm{j}}$ for the single patches, and will further determine the global collocation points in the multi-patch isogeometric collocation technique in Section~\ref{subsec:biharmonic}. 

\begin{figure}
    \centering
    \includegraphics[scale=0.43]{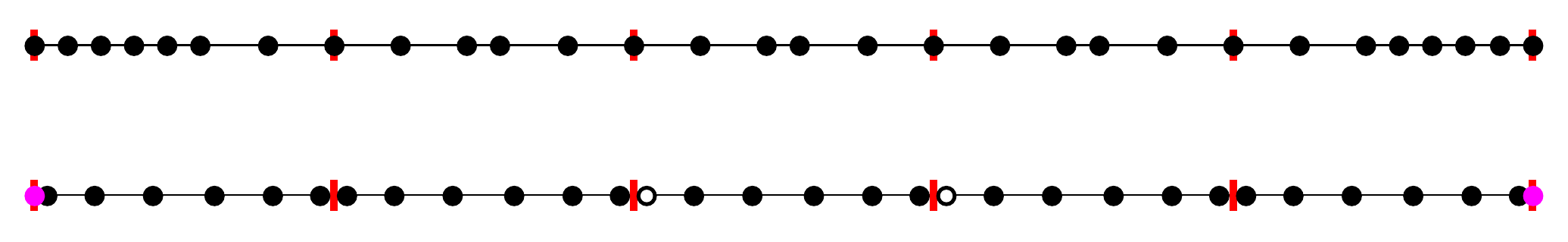}
    \caption{Collocation points for the case of $p=9$, $r=4$ and $k=4$. Top row: Greville points. Bottom row: All superconvergent points (black and white points) and clustered superconvergent points (only black points). The violet points are the two boundary points which we have to add to the set of superconvergent points. 
    {\new The vertical red ticks represent the knots.}
    }
    \label{fig:GrevilleSuperconvergent9}
\end{figure}

\paragraph{The Greville points} The first choice of points are the Greville points which are defined for an arbitrary $p$ and $r$ as 
\[
 \loc_{j}^{p,r}=\frac{t_{j+1}^{p,r} + \ldots t_{j+p}^{p,r}}{p}, \quad  j \in \{0,1,\ldots, n-1 \}, 
\]
and which are shown for the case $p=9$, $r=4$ and $k=4$ in Fig.~\ref{fig:GrevilleSuperconvergent9} (top row). 

\paragraph{The superconvergent points} The second choice of points are the so-called superconvergent points which are defined in our case of the biharmonic equation as the roots of its Galerkin residual $D^4(u-u_h)$. Estimates of these points can be computed by solving the simple biharmonic equation  
\begin{equation} \label{eq:biharmonic_equation}
u^{(4)}(x) =f(x), \quad x \in (0,1), \quad {\rm and} \quad u(0)=u(1)=0, \quad u'(0)=u'(1)=0, 
\end{equation}
with some particular function~$f$ as demonstrated in~\cite{GomezLorenzisVariationalCollocation}. Superconvergent points for fourth order problems have been studied so far just for the one-patch case~\cite{GomezLorenzisVariationalCollocation,RealiGomez2015,Maurin2018}, and there only for splines with maximal regularity, i.e. $r=p-1$. However, in the multi-patch case the use of splines with maximal regularity is not possible, and splines with a regularity $4 \leq r \leq p-5$ have to be employed instead. Superconvergent points for splines with a lower regularity than $r=p-1$ have been investigated so far only for the case of second order problems~\cite{KaVi20}, 
and will be extended below to fourth order problems for the case of the underlying spline space~$\mathcal{S}^{9,4}_h([0,1])$. 

Following the approach~\cite{GomezLorenzisVariationalCollocation} of studying the Galerkin residual $D^4(u-u_h)$ of a particular biharmonic problem~\eqref{eq:biharmonic_equation}, we get for the spline space~$\mathcal{S}^{9,4}_h([0,1])$ the superconvergent points on each knot span with respect to the reference interval~$[-1,1]$ as the roots of the sextic polynomial
$$
4823 x^6 - 5915 x^4 + 1665 x^2 - 61.
$$
Numerical approximations of the six roots are {\new given in Tab.~\ref{tab:superconvergent points} (first row).} 
All univariate superconvergent points are presented in Fig.~\ref{fig:GrevilleSuperconvergent9} (bottom row). 
Since the boundary points of the domain interval are not contained in the set of superconvergent points, they have to be added to the set to be able to impose the required boundary conditions.
To get the same number of collocation points as the dimension of the spline space, it can be necessary 
to select a proper subset of superconvergent points, or conversely, 
to add some points 
to the set of superconvergent points, cf.~\cite{MonSanTam2017,KaVi20}.
In our case the difference between the number of all superconvergent points and the dimension of the underlying spline space $\mathcal{S}_{h}^{9,4}([0,1])$ is 
$$
  6(k+1)+2 - \dim \mathcal{S}_{h}^{9,4}([0,1]) = 6(k+1)+2 - (5k+10) = k-2,
$$
so we have to omit $k-2$ points.
For this purpose, we skip for all $k-2$ most inner knots (at least two away from the boundary) the closest collocation point to the right of the knot, compare Fig.~\ref{fig:GrevilleSuperconvergent9} (white points). For special cases $k=0,1$, we either add the inner knot (case $k=1$) or we add the second and the second last Greville point (case $k=0$). 

\paragraph{Separation of collocation points} Above, we have explained for the spline space~$\mathcal{S}_{h}^{9,4}([0,1])$ two possible choices of univariate local collocation points~$\loc_j$, namely the Greville and superconvergent points, where in both cases the number of points coincides with the dimension of the spline space. Recall that the univariate local collocation points~$\loc_j$ define via $\bfm{\loc}_{\bfm{j}} = (\loc_{j_1},\loc_{j_2})$ the local collocation points $\bfm{\loc}_{\bfm{j}} \in [0,1]^2$. Then, each local collocation point $\bfm{\loc}_{\bfm{j}}$ for a patch~$\overline{\Omega^{(i)}}$ defines via $\bfm{\glob}_{i_{\bfm{j}}} = \ab{F}^{(i)}(\bfm{\loc}_{\bfm{j}})$ the global collocation point $\bfm{\glob}_{i_{\bfm{j}}} \in \overline{\Omega^{(i)}}$, $i_{\bfm{j}} \in \mathcal{J}$, taking under consideration, that if the global collocation point would lie on more than one patch, the global collocation point is just selected once, namely for the patch~$\overline{\Omega^{(i)}}$ with the lowest patch index~$i$, cf.~\eqref{eq:patch_index_selection}. 

\begin{figure}
    \centering
    \begin{minipage}{7cm}
        \includegraphics[scale=0.39]{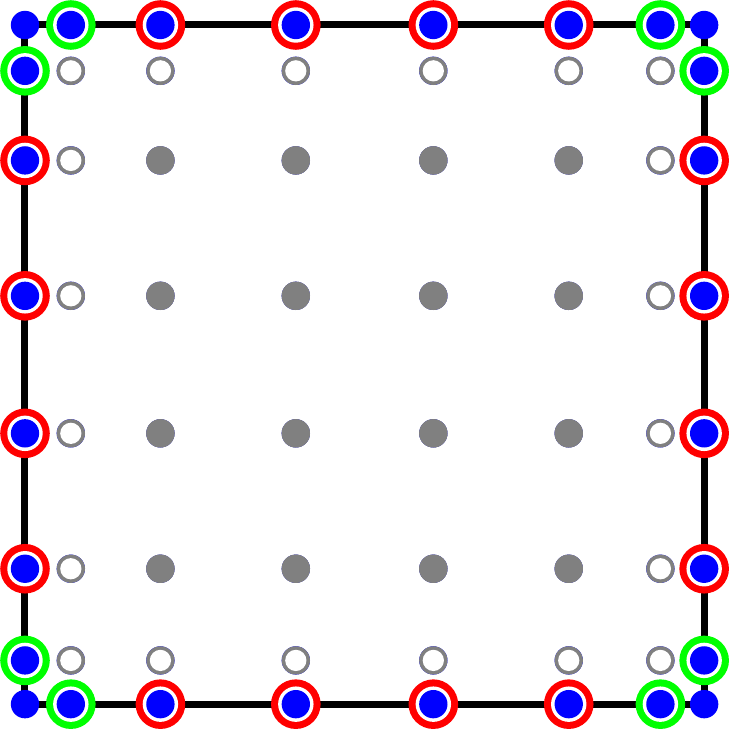}
    \end{minipage}
     \begin{minipage}{7cm}
     \includegraphics[scale=0.21]{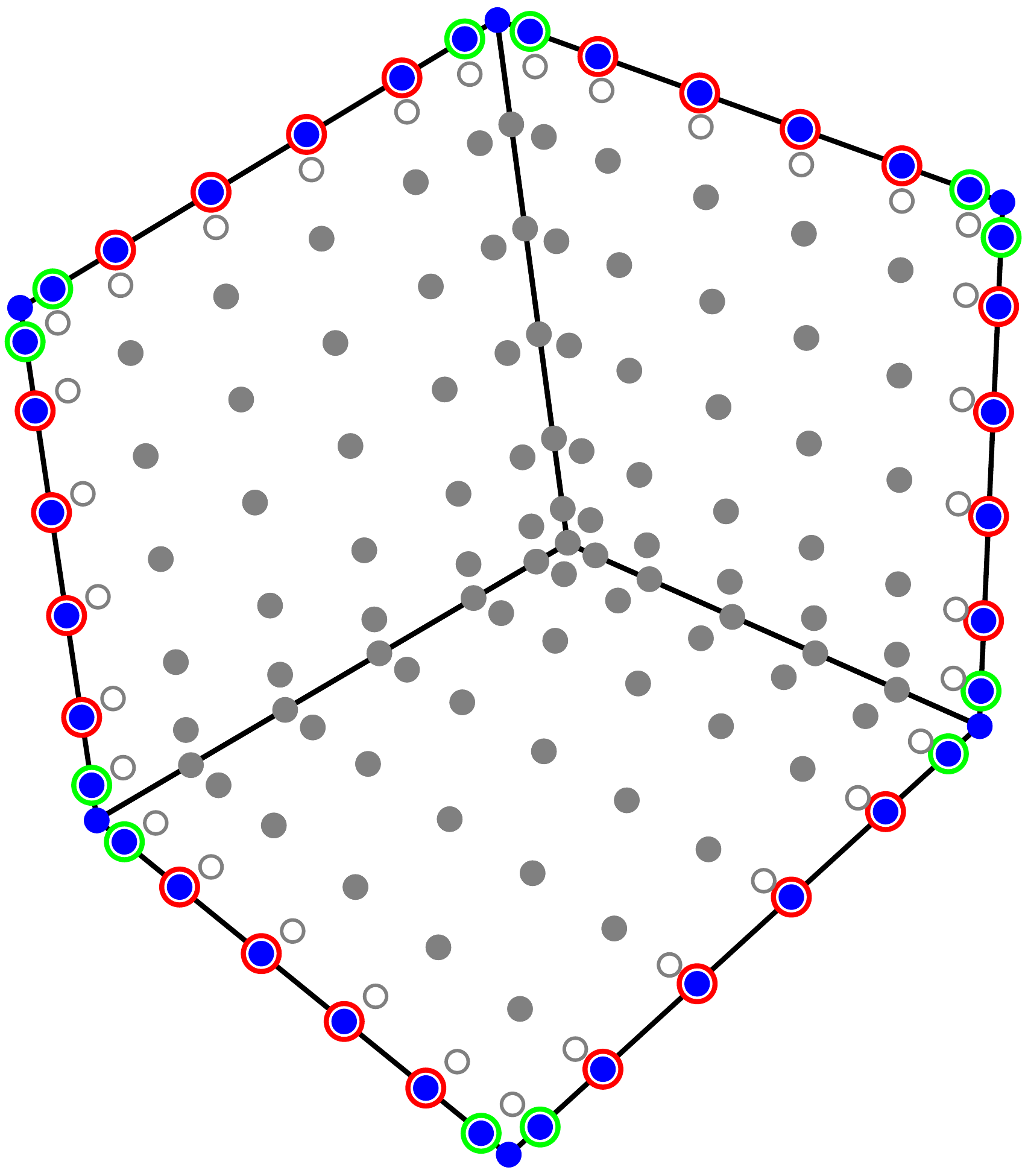}
     \end{minipage}
    \caption{Separation of collocation points into inner and boundary collocation points. Gray inner collocation points are the ones where the biharmonic equation is imposed. The empty inner collocation points are the ones that are omitted. Blue boundary collocation points are points where the boundary values are interpolated, while the interpolation of normal derivatives is related to the red and green circles, where the green ones denote points where the two equations are averaged. The left figure represents the situation for the one-patch case, while the right one shows the generalization to the multi-patch case.}
    \label{fig:collocationPointsOnePatch}
\end{figure}

The usage of the multi-patch isogeometric collocation approach presented in Section~\ref{subsec:biharmonic} still requires the separation of the set of global collocation points~$\bfm{\glob}_j$, $j \in \mathcal{J}$, into inner collocation points $\bfm{\glob}_j^I$, $j \in \mathcal{J}_I$, and boundary collocation points $\bfm{\glob}_j^B$,  $j \in \mathcal{J}_B$. The straightforward separation would be to take all global collocation points which lie on the boundary of the domain $\partial \Omega$ as boundary collocation points, 
and to take the remaining points as inner collocation points. 
However, this separation would lead even in the one-patch case to an overdetermined linear system, since we would have too many collocation points in total. Therefore, we will follow the separation technique used in~\cite{RealiGomez2015, GomezRealiSangali2014} for the one-patch case, which will lead in this case to a square linear system, and will extend it to the multi-patch case, where we will obtain just a slightly overdetermined linear system, which will then be solved by means of the least-squares approach in the numerical examples in Section~\ref{section_Numerical_examples}. The separation of the collocation points will work as follows. We will omit some of the inner collocation points where the biharmonic equation is imposed, and will choose the boundary collocation points for interpolation of the normal derivatives in a proper way. In the one-patch case, cf. Fig.~\ref{fig:collocationPointsOnePatch} (left), we will omit the last outer ring of inner collocation points, see~Fig.~\ref{fig:collocationPointsOnePatch} (white points), and we will further omit eight interpolation conditions at the boundary. For the later purpose, 
we will not interpolate the normal derivatives at the four corner points, and will additionally
average the normal derivative equations corresponding to the two boundary points closest to the corner points, see Fig.~\ref{fig:collocationPointsOnePatch} (green circles).  This strategy will be generalized in a straightforward way to the multi-patch case as demonstrated in Fig.~\ref{fig:collocationPointsOnePatch} (right). 

\paragraph{The convergence behavior of the collocation points} For fourth order problems, the convergence behavior of both types of collocation points has been studied so far only with respect to the $L^2$-norm for the one-patch case and there just using splines with maximal regularity, where for odd spline degrees~$p$ a convergence of order $\mathcal{O}(h^{p-3})$ and $\mathcal{O}(h^{p-2})$ for the Greville and superconvergent points, respectively, have been reported, cf.~\cite{GomezLorenzisVariationalCollocation,Reali2015}. In the following section, we will first numerically demonstrate for the one-patch case that the convergence behavior of both types of collocation points will be unchanged by employing the spline space~$\mathcal{S}^{p,r}_h([0,1])$ with $p=9$ and $r=4$. Moreover, we will study for this spline space the convergence behavior with respect to (equivalents of) the $H^s$-seminorms for $1 \leq s \leq 4$, where in case of the Greville points a convergence of order~$\mathcal{O}(h^{p-3})$ for all seminorms will be observed, and where in case of the superconvergent points the convergence will be of order $\mathcal{O}(h^{p-2})$ except for the equivalent of the $H^4$-seminorm with a convergence order~$\mathcal{O}(h^{p-3})$. In the multi-patch case, we will present several examples which will indicate for both considered types of collocation points, i.e. the Greville and superconvergent points, that the convergence rates for all (semi)norms will be the same as in the one-patch case.

\section{Numerical examples}  \label{section_Numerical_examples}

We will demonstrate on the basis of several examples the potential of our isogeometric collocation method for solving the biharmonic equation over planar multi-patch domains, and will study the convergence behavior under $h$-refinement for both types of considered collocation points, namely for the Greville and superconvergent points. While the first example (Example~\ref{ex:onepatchdomain}) will deal with the one-patch case, the remaining examples (Examples~\ref{ex:threepatchdomain}--\ref{ex:twopatchdomainSmoothness3}) will study different instances of multi-patch domains. As already mentioned in the previous section, the use of both sets of collocation points will lead in case of a multi-patch domain (Examples~\ref{ex:threepatchdomain}--\ref{ex:twopatchdomainSmoothness3}) to a slightly overdetermined system of linear equations, {\new see Tab.~\ref{tab:matrixSizes} for the dimensions of these resulting collocation matrices.}
{\new This overdetermined linear system could be solved e.g.
by using diagonal scaling (\cite{Br95}) on a symmetric positive definite matrix obtained by employing the normal system approach or by using QR-decomposition. Condition numbers with estimated growth rates 
for both approaches and for all considered examples are presented in Fig.~\ref{fig:conditionNumbers}. Since the QR decomposition approach gives lower condition numbers, we will use this method to construct all multi-patch examples below.} 

In all examples (Examples~\ref{ex:onepatchdomain}--\ref{ex:twopatchdomainSmoothness3}), we will solve the biharmonic equation~\eqref{eq:biharmonic}, where the right side functions $g, g_1$ and $g_2$ are obtained from the exact solution 
\begin{equation}  \label{eq:exactSolution}
 u(x_1,x_2)= \cos\left(\frac{x_1}{2}\right) \sin\left(\frac{x_2}{2}\right). 
\end{equation}
In the first three examples (Examples~\ref{ex:onepatchdomain}--\ref{ex:fourpatchdomainHole}), we will employ the $C^4$-smooth isogeometric spline spaces~$\W^4$ with a spline degree~$p=9$, with an inner patch regularity~$r=4$ and with mesh sizes $h=1/2^i$, $i=2,\ldots,5$, as discretization spaces to compute a $C^4$-smooth approximation of the solution of the biharmonic equation~\eqref{eq:biharmonic}, and will denote the resulting space $\W^4$ for a specific mesh size $h$ by $\mathcal{W}^4_h$. We will investigate the quality of the obtained approximant $u_h \in \mathcal{W}^4_h$ with respect to the exact solution \eqref{eq:exactSolution} by considering the relative errors with respect to the $L^2$-norm and $H^1$-seminorm as well as 
with respect to equivalents (cf. \cite{BaDe15})
of the $H^s$-seminorms, $2\leq s \leq 4$, i.e.
\begin{equation} \label{eq:eqiuv2seminorms}
\frac{\| u-u_h\|_{L^2}}{\| u \|_{L^2}}, \,
\frac{| u-u_h|_{H^1}}{| u |_{H^1}}, \,
 \frac{\| \Delta u- \Delta u_h\|_{L^2}}{\| \Delta u \|_{L^2}}, \,
\frac{ \| \nabla (\Delta u) - \nabla(\Delta u_h)\|_{L^2}}{\| \nabla(\Delta u )\|_{L^2}}, \,
\frac{\| \Delta^2 u- \Delta^2 u_h\|_{L^2}}{\| \Delta^2 u \|_{L^2}}.
\end{equation}
For the sake of brevity, we will refer in this section to the relative errors with respect to the equivalents of the $H^s$-seminorms, $2\leq s \leq 4$, just as relative errors with respect to the corresponding $H^s$-seminorms. 

\begin{ex} \label{ex:onepatchdomain}
We consider the bilinear one-patch domain $\overline{\Omega}$ shown in Fig.~\ref{fig:onePatchErrors} (top left). In this case, the $C^4$-smooth discretization spaces~$\W_h^4$, $h=1/2^i$, $i=2,\ldots,5$, technically consist of only 
patch functions (cf. last paragraph in Section~\ref{C4_isogeometric_splinespace}), and the resulting linear system~\eqref{eq:collocationSystemLocal} is a square system.
The resulting relative errors~\eqref{eq:eqiuv2seminorms} 
by using the Greville and superconvergent points are shown in Fig.~\ref{fig:onePatchErrors}.

\begin{figure}[hbt!]
    \centering
    \includegraphics[scale=0.20]{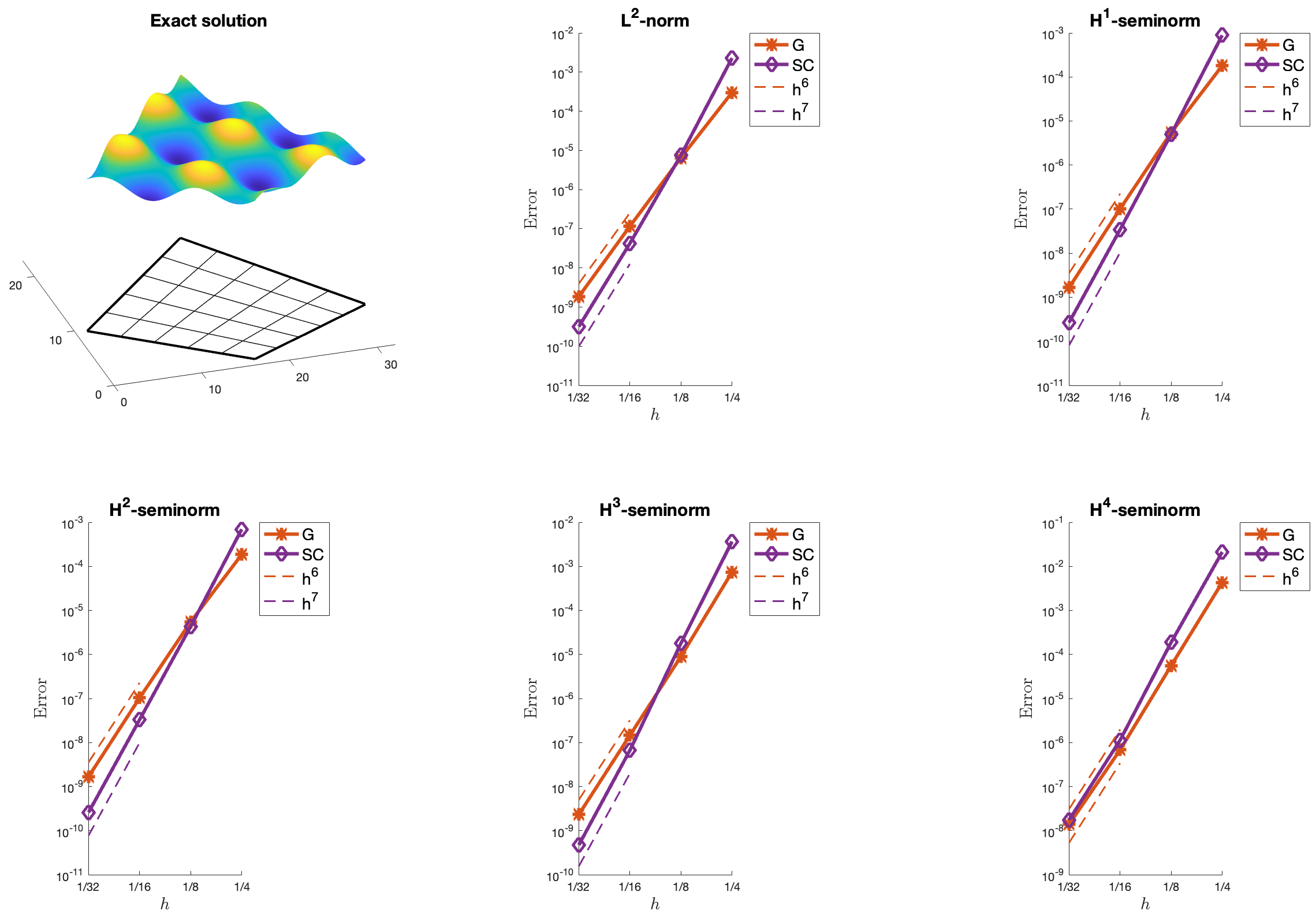}    
    \caption{Example~\ref{ex:onepatchdomain}. Plots of the 
    bilinear one-patch domain, of the considered exact solution~\eqref{eq:exactSolution} as well as of the relative errors~\eqref{eq:eqiuv2seminorms} computed with respect to the $L^2$-norm and with respect to 
    the $H^s$-seminorms, $1\leq s \leq 4$, {\new using the Greville (G) and the superconvergent points (SC) as collocation points. 
    } }
    \label{fig:onePatchErrors}
\end{figure}
The convergence rate for the relative $L^2$-error coincides with the known results for splines with maximal regularity, i.e. $r=p-1$, for the one-patch  case~\cite{GomezLorenzisVariationalCollocation,Reali2015}, and is in case of odd spline degree $p=9$ of order~$\mathcal{O}(h^{p-3})=\mathcal{O}(h^6)$ for the Greville points and of order $\mathcal{O}(h^{p-2})=\mathcal{O}(h^7)$ for the superconvergent points. Furthermore, the error plots indicate the convergence orders with respect  
to the $H^s$-seminorms, $1 \leq s \leq 4$, where for the Greville points the order remains $\mathcal{O}(h^{p-3})=\mathcal{O}(h^6)$ for all seminorms, while in the case of the superconvergent 
points the orders are $\mathcal{O}(h^{p-2})=\mathcal{O}(h^7)$ for all the seminorms except for the $H^4$-seminorm, where $\mathcal{O}(h^{p-3})=\mathcal{O}(h^6)$ is already the best possible one, which would be also obtained in  
case of the Galerkin approach. 
\end{ex}

\begin{ex} \label{ex:threepatchdomain}
In this example we consider the bilinear three-patch and five patch domains~$\overline{\Omega}$ visualized in Fig.~\ref{fig:threePatchErrors} (top left) and Fig.~\ref{fig:fivePatchErrors} (top left), respectively.
\begin{figure}[hbt!]
    \centering
    \includegraphics[scale=0.2]{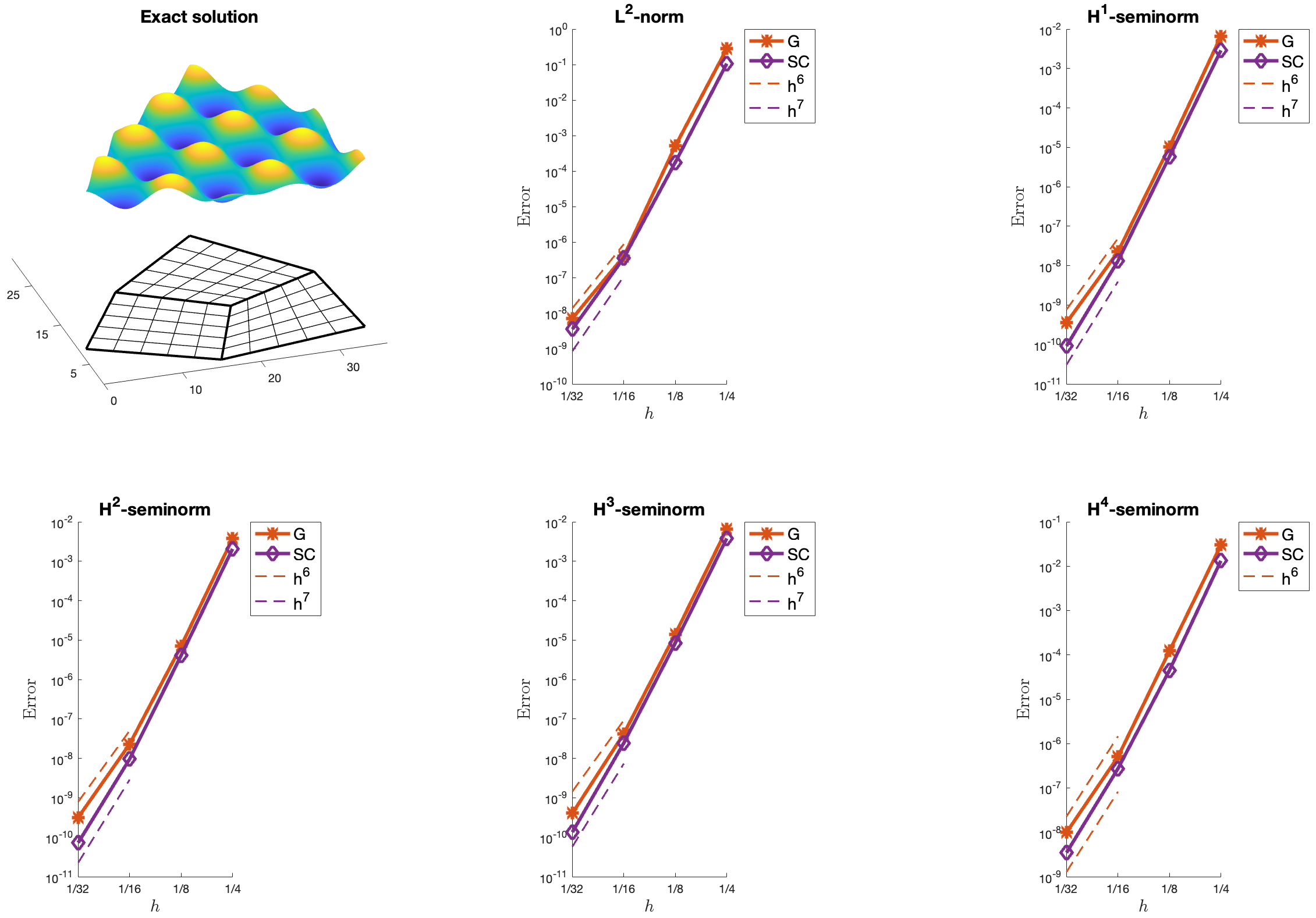}
    \caption{Example~\ref{ex:threepatchdomain} (Three-patch case). Plots of the bilinear three-patch domain, of the considered exact solution~\eqref{eq:exactSolution} as well as of the relative errors~\eqref{eq:eqiuv2seminorms} computed with respect to the $L^2$-norm and with respect to the $H^s$-seminorms, $1\leq s \leq 4$, {\new using the Greville (G) and the superconvergent points (SC) as collocation points. 
    } }
    \label{fig:threePatchErrors}
\end{figure}
The $C^4$-smooth discretization spaces~$\mathcal{W}^4_h$, $h=1/2^i$, $i=2,\ldots,5$, are constructed as presented in Section~\ref{C4_isogeometric_splinespace} following the mentioned adaptation from the last paragraph for the mesh sizes $h=\frac{1}{4}$ and $h=\frac{1}{8}$.  
Again, we compare the resulting relative errors~\eqref{eq:eqiuv2seminorms} with respect to the $L^2$-norm and 
to $H^s$-seminorms, $1 \leq s \leq 4$,
by using the Greville and superconvergent points, see Fig.~\ref{fig:threePatchErrors} for the three-patch domain and Fig.~\ref{fig:fivePatchErrors} for the five-patch domain.
The numerical results indicate in case of the Greville points a convergence of order~$\mathcal{O}(h^{p-3}) = \mathcal{O}(h^{6}) $ independent of the considered (semi)norm, and show in case of the superconvergent points an improved convergence of order~$\mathcal{O}(h^{p-2}) = \mathcal{O}(h^{7})$ for all (semi)norms except for 
the $H^4$-seminorm, where the order $\mathcal{O}(h^{p-3})=\mathcal{O}(h^{6})$ is already optimal,  
{\new 
which is for all (semi)norms the same convergence behavior as in the one-patch case in Example~\ref{ex:onepatchdomain}.
}

\begin{figure}[hbt]
    \centering
    \includegraphics[scale=0.2]{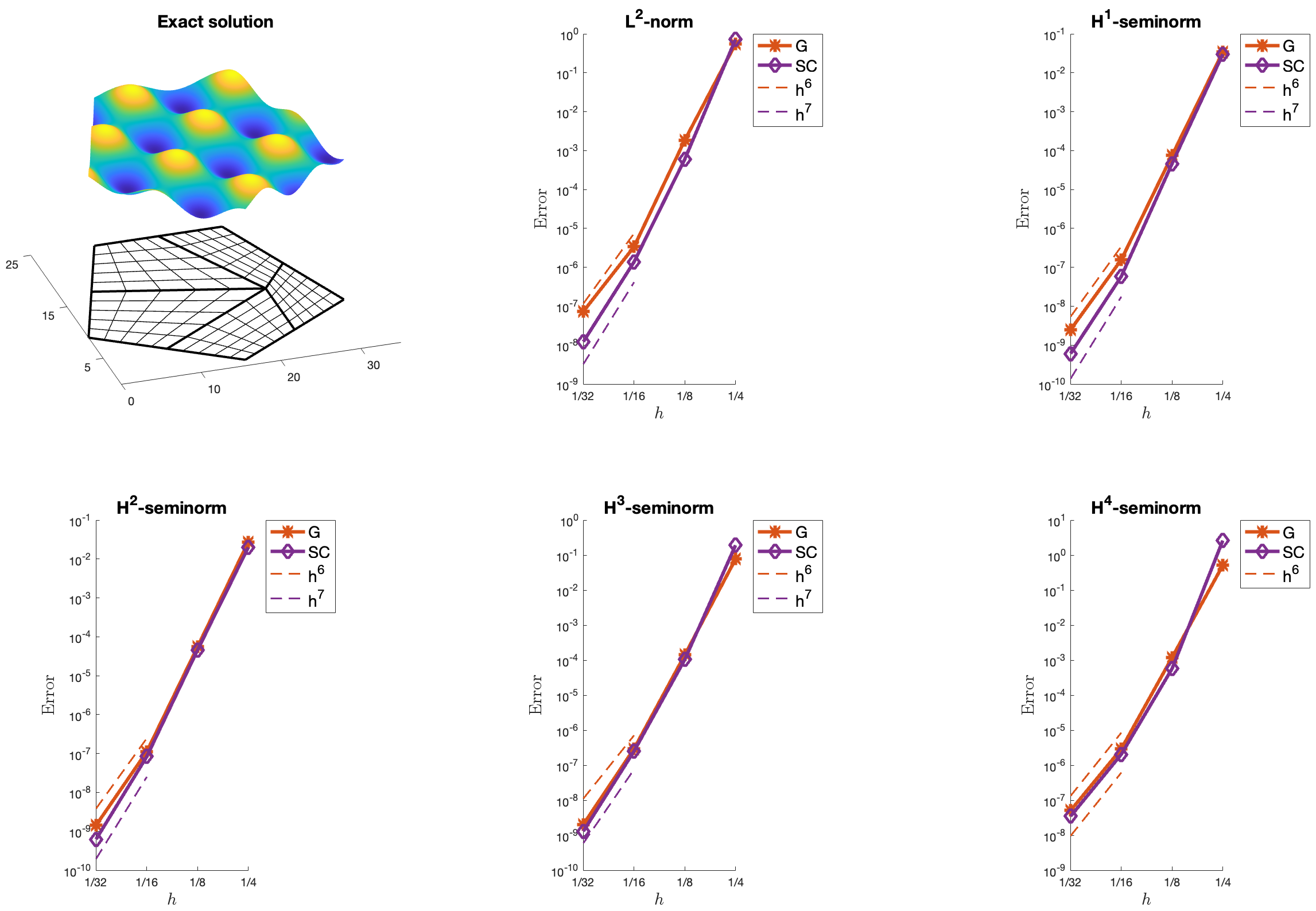}
    \caption{Example~\ref{ex:threepatchdomain} (Five-patch case). Plots of the five-patch domain, of the considered exact solution~\eqref{eq:exactSolution} as well as of the relative errors~\eqref{eq:eqiuv2seminorms} computed with respect to the $L^2$-norm and with respect to the $H^s$-seminorms, $1\leq s \leq 4$, {\new using the Greville (G) and the superconvergent points (SC) as collocation points. 
    } }
    \label{fig:fivePatchErrors}
\end{figure}
\end{ex}

\begin{ex} \label{ex:fourpatchdomainHole}
{\new
We consider now a generalization of bilinear multi-patch domains to the case of bilinear-like $G^4$ multi-patch geometries, see e.g.~\cite{CoSaTa16, KaVi17c, KaVi20b}. These multi-patch parameterizations possess as in the bilinear case linear gluing functions and also allow the construction of globally $C^4$-smooth isogeometric spline spaces with optimal approximation properties. 
More precisely, a multi-patch parameterization~$\ab{F}$ is called \emph{bilinear-like $G^{4}$} if for any two neighboring 
patches~$\Omega^{(i_0)}$ and $\Omega^{(i_1)}$, $i_0,i_1 \in \mathcal{I}_{\Omega}$, with~$\overline{\Gamma^{(i)}} = \overline{\Omega^{(i_0)}} \cap 
\overline{\Omega^{(i_1)}}$, $i \in \mathcal{I}_{\Gamma}^I$, 
and corresponding geometry mappings $\ab{F}^{(i_0)}$ and $\ab{F}^{(i_1)}$, $i_0,i_1 \in \mathcal{I}_{\Omega}$, parameterized as in Fig.~\ref{fig:twopatchCase},
there exist linear functions 
$\alpha^{(i,\Side)},\beta^{(i,\Side)}, \Side \in \{\LL,\RR \}$, such that
 \begin{equation*}   \label{eq:FC}
 \ab{F}_\ell^{(\LL)}(\xi) = \ab{F}_\ell^{(\RR)}(\xi) =: \ab{F}_\ell(\xi), \quad \ell =0,1,\ldots,\sm,
 \end{equation*}
 with
 \begin{equation*}   \label{eq:FC2}
 \ab{F}_\ell^{(\Side)}(\xi) = \left(\alpha^{(i,\Side)}(\xi)\right)^{-\ell}\, \partial_1^\ell \ab{F}^{(\Side)}(0,\xi) - \sum_{i=0}^{\ell-1} {\ell \choose i} 
 \left(\frac{\beta^{(i,\Side)}(\xi)}{\alpha^{(i,\Side)}(\xi)}\right)^{\ell-i}  \dd^{\ell-i} \ab{F}_i^{}(\xi) ,\quad \Side\in \{\LL,\RR\}.
 \end{equation*}
The advantage of using bilinear-like $G^4$ multi-patch parameterizations instead of just bilinear multi-patch parameterizations is the possibility to deal with multi-patch domains having curved interfaces and 
boundaries, see e.g.~\cite{KaSaTa17b,KaSaTa19b,KaSaTa19a, KaVi20b}. 

We study more precisely the 
{\new bilinear-like $G^4$} four-patch domain~$\overline{\Omega}$ shown in Fig.~\ref{fig:fourPatchErrors} (top left), that consists of polynomial patches of bi-degree~$9$, and possesses a $C^4$-smooth outer boundary and an inner boundary with sharp corners. The domain can represent a real engineering domain, namely a washer with a square hole. Compared to round washers, square washers have a higher surface area, which improves torque distribution. They help in corrosion resistance and rotation prevention.} 
The $C^4$-smooth discretization spaces~$\mathcal{W}^4_h$, 
$h=1/2^i$, $i=2,\ldots,5$ technically consist now of only 
patch and edge functions (cf.~last paragraph in Section~\ref{C4_isogeometric_splinespace}).~We compare the resulting relative errors~\eqref{eq:eqiuv2seminorms} with respect to the $L^2$-norm and with respect to 
the $H^s$-seminorms, $1 \leq s \leq 4$, 
by using the Greville and superconvergent points as 
collocation points, see Fig.~\ref{fig:fourPatchErrors}. {\new Again}, the numerical results indicate the same convergence behavior as in the 
{\new previous examples.}
\begin{figure}[htb!]
    \centering
    \includegraphics[scale=0.2]{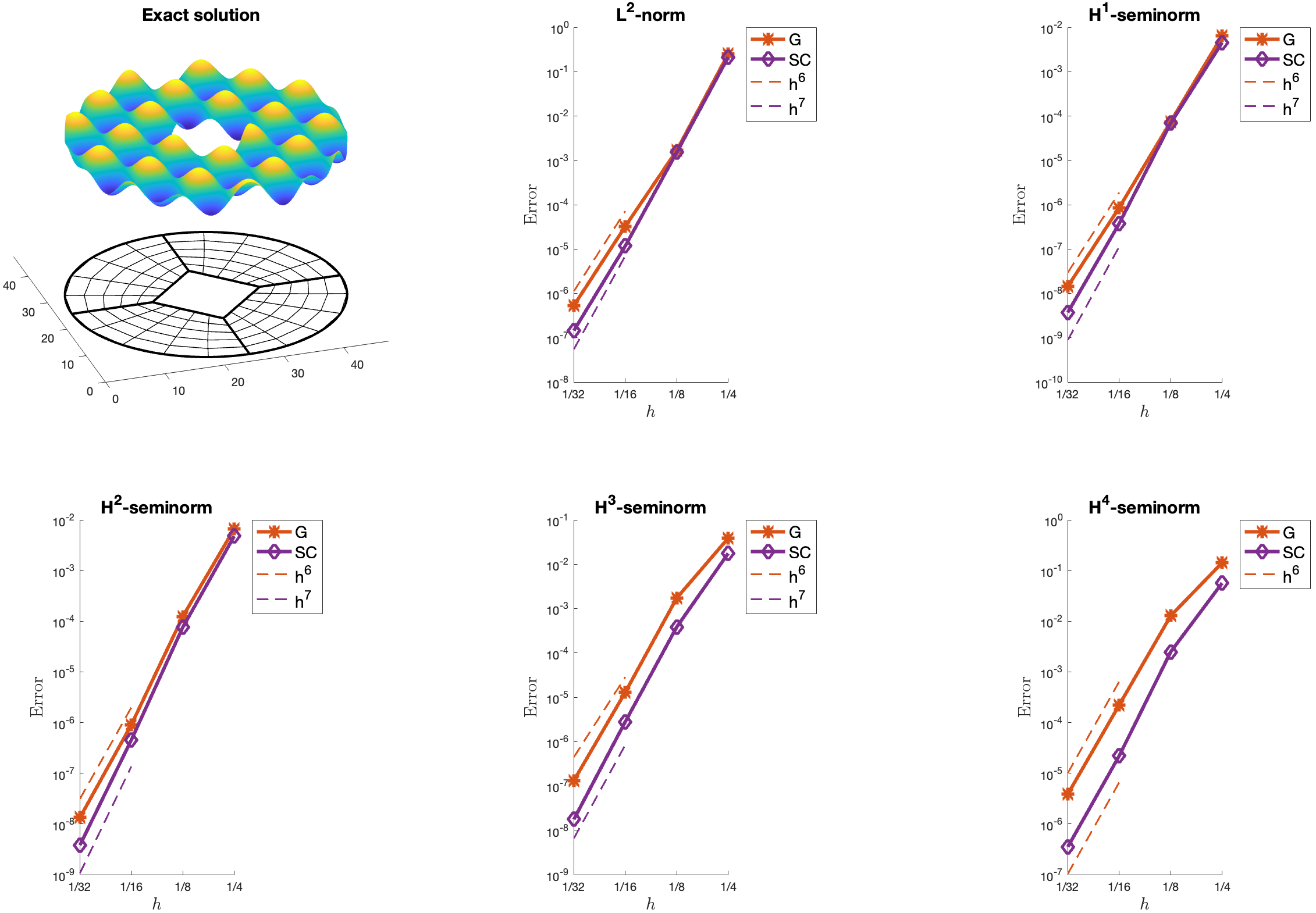}    
    \caption{Example~\ref{ex:fourpatchdomainHole}. Plots of the {\new bilinear-like $G^4$} four-patch domain, of the considered exact solution~\eqref{eq:exactSolution} as well as of the relative errors~\eqref{eq:eqiuv2seminorms} computed with respect to the $L^2$-norm and with respect to the $H^s$-seminorms, $1\leq s \leq 4$, {\new using the Greville (G) and the superconvergent points (SC) as collocation points. 
    } }
    \label{fig:fourPatchErrors}
\end{figure}
\end{ex}

In the last example, we will present a first possible way to slightly reduce the required spline degree~$p$ of the used discretization space by reducing the regularity of the numerical solution as well, cf.~\cite{SuperConvergent2015}. For this purpose, we will relax the condition of computing a $C^4$-smooth approximation of the solution of the biharmonic equation by constructing instead just a $C^3$-smooth approximation. Since the obtained approximant will now be globally only $C^3$-smooth, we will have to choose the collocation points carefully in order to avoid collocating at the knots and common edges~$\overline{\Gamma^{(i)}}$, $i \in \mathcal{I}_{\Gamma}^{I}$, where the fourth derivative of the approximant may be discontinuous. Again, we will study the convergence rates of the relative errors~\eqref{eq:eqiuv2seminorms}.

\begin{ex} \label{ex:twopatchdomainSmoothness3}
Let $\overline{\Omega}$ be the bilinear two-patch {\new L-shape} domain given in Fig.~\ref{fig:twoPatchErrors} (top left). The goal is to compute via our collocation approach a $C^3$-smooth approximant $u_h \in \W_h^3$ of the solution of the biharmonic equation, where $\W^3_h$ is the $C^3$-smooth isogeometric discretization space~\cite{KaVi20b} for a specific mesh size $h=1/2^i$, $i=2,\ldots,5$, whose construction is similar to the one for the $C^4$-smooth spline space~$\W^4_h$ described in Section~\ref{C4_isogeometric_splinespace}. The $C^3$-smooth spline space~$\W^3_h$
allows the use of underlying spaces $\mathcal{S}_{h}^{p,3}([0,1])$ with $p \geq 7$.  
In this example, we will consider underlying spline spaces $\mathcal{S}_{h}^{p,3}([0,1])$ with $p=7$ and $p=8$ for both sets of collocation points, namely for the Greville and superconvergent points.
Since the regularity~$r=3$ is odd and less than $\frac{p}{2}$, it is easy to see that the Greville points for the two spline spaces $\mathcal{S}_{h}^{7,3}([0,1])$ and $\mathcal{S}_{h}^{8,3}([0,1])$ never coincide with inner knots, cf.~Fig.~\ref{fig:GrevilleSuperconvergent8and7}. 
\begin{figure}[htb!]
    \centering
    \includegraphics[scale=0.43]{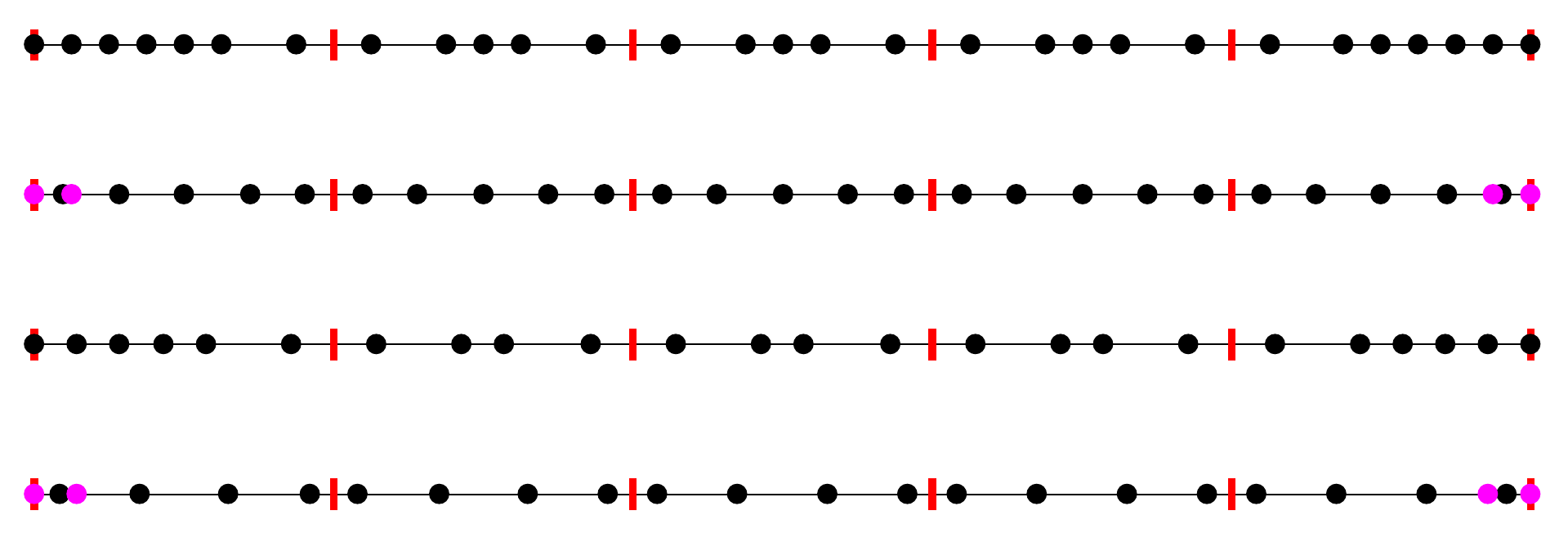}
    \caption{Greville points (first and third row) and all superconvergent points (second and fourth row). The first two rows correspond to the case $(p,r)=(8,3)$, while the last two rows to the case $(p,r)=(7,3)$. The violet points are the two pairs of boundary Greville points which we have to add to the set of superconvergent points. {\new The vertical red ticks represent the knots.}
    }
    \label{fig:GrevilleSuperconvergent8and7}
\end{figure}
Therefore, we only have to avoid collocating at the common edge of the two-patch domain, where the solution is not $C^4$-smooth. Similarly as for the case $(p,r)=(9,4)$ in Section~\ref{sec:SelectionCollocPoints}, we can construct 
the superconvergent 
points for the cases $(p,r)=(8,3)$ and $(p,r)=(7,3)$. In the first case, i.e. for the spline space~$\mathcal{S}^{8,3}_h([0,1])$, the superconvergent points 
on each knot segment (with respect to the reference interval~$[-1,1]$) are the roots of the quintic polynomial 
$$
3003 x^5 - 2530 x^3 + 375 x,
$$
while in the second case, i.e. for the spline space~$\mathcal{S}^{7,3}_h([0,1])$, 
we have to take the roots of the quartic polynomial 
$$
495 x^4 - 390 x^2 + 31,
$$
{\new cf.~Tab.~\ref{tab:superconvergent points}.}
\begin{table}[htb!] 
    \centering
    {\new 
    \begin{tabular}{c|c} 
    \hline\\[-0.3cm]
       $p=9$, $r=4$  & $\pm 0.2072795685478027, \; \pm 0.5963052503103114, \; \pm 0.9098737952346008$ \\
        \hline\\[-0.3cm]      
        $p=8$, $r=3$  &  $0, \; \pm \sqrt{5 (253 - 2 \sqrt{4741})/3003}, \; \pm \sqrt{5 (253 + 2 \sqrt{4741})/3003} $\\
         \hline\\[-0.3cm]  
       $p=7$, $r=3$  &  $\pm \sqrt{(65 - 6 \sqrt{70})/165}, \; \pm \sqrt{(65 + 6 \sqrt{70})/165} $\\
       \hline
    \end{tabular}
    \caption{Sets of superconvergent points on each knot segment with respect to the reference interval~$[-1,1]$ for all three cases considered in the paper.}
    \label{tab:superconvergent points}
    }
\end{table}

Again, we have to find a proper set of superconvergent points whose cardinality coincides with the dimension of the underlying spline space~$\mathcal{S}^{8,3}_h([0,1])$ or $\mathcal{S}^{7,3}_h([0,1])$. In this example, we have to add (in both cases) the first two Greville points at both ends of the domain interval, see Fig.~\ref{fig:GrevilleSuperconvergent8and7}. The same strategy has been already used for selecting additional points close to the boundary in \cite{GomezLorenzisVariationalCollocation, KaVi20, MonSanTam2017}. 
When considering our two-patch domain (Fig.~\ref{fig:twoPatchErrors}, top left), we only have, as for the Greville collocation points, to avoid collocating at the common edge,
where the numerical solution is not $C^4$-smooth. 
\begin{figure}[htb!]
    \centering
    \includegraphics[scale=0.2]{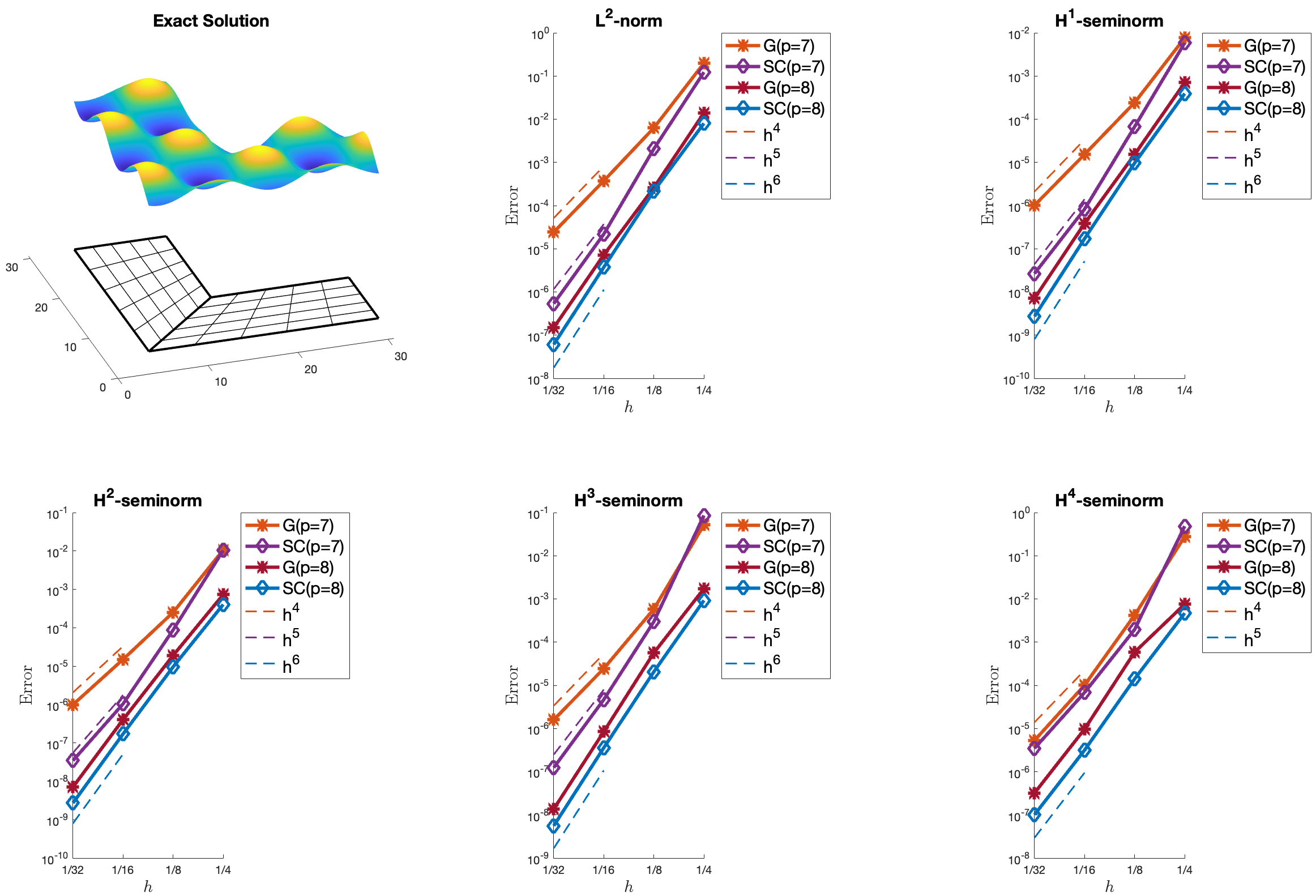}  
    \caption{Example~\ref{ex:twopatchdomainSmoothness3}. Plots of the bilinear two-patch domain, of the considered exact solution~\eqref{eq:exactSolution} as well as of the relative errors~\eqref{eq:eqiuv2seminorms} computed with respect to the $L^2$-norm and with respect to the $H^s$-seminorms, $1\leq s \leq 4$, {\new using the Greville (G) and the superconvergent points (SC) as collocation points. 
    } }
    \label{fig:twoPatchErrors}
\end{figure}

The convergence rates for the relative errors~\eqref{eq:eqiuv2seminorms} with respect to the $L^2$-norm and 
with respect to the $H^s$-seminorms, $1 \leq s \leq 4$, 
are presented in Fig.~\ref{fig:twoPatchErrors}. For odd spline degree~$p=7$, the orders of convergence coincide with the ones in the previous examples, i.e. they are of order~$\mathcal{O}(h^{p-3})$ in the case of Greville points, and they are of order~$\mathcal{O}(h^{p-2})$ for all (semi)norms except for 
the $H^4$-seminorm, where the order is $\mathcal{O}(h^{p-3})$, in the case of superconvergent points. For even spline degree~$p=8$, the orders of convergence are the same for the Greville and superconvergent points, namely $\mathcal{O}(h^{p-2})$ for the $L^2$-norm and for the $H^s$-seminorms, for $1\leq s \leq 3$, and $\mathcal{O}(h^{p-3})$ for the $H^4$-seminorm. Note that the observed convergence behavior with respect to the $L^2$-norm is for even spline degree again in agreement with the known one-patch result for splines with maximal regularity, cf.~\cite{GomezLorenzisVariationalCollocation,Reali2015}.
\end{ex}

{\new 
\begin{figure}[htb!]
    \centering
    \includegraphics[scale=0.188]{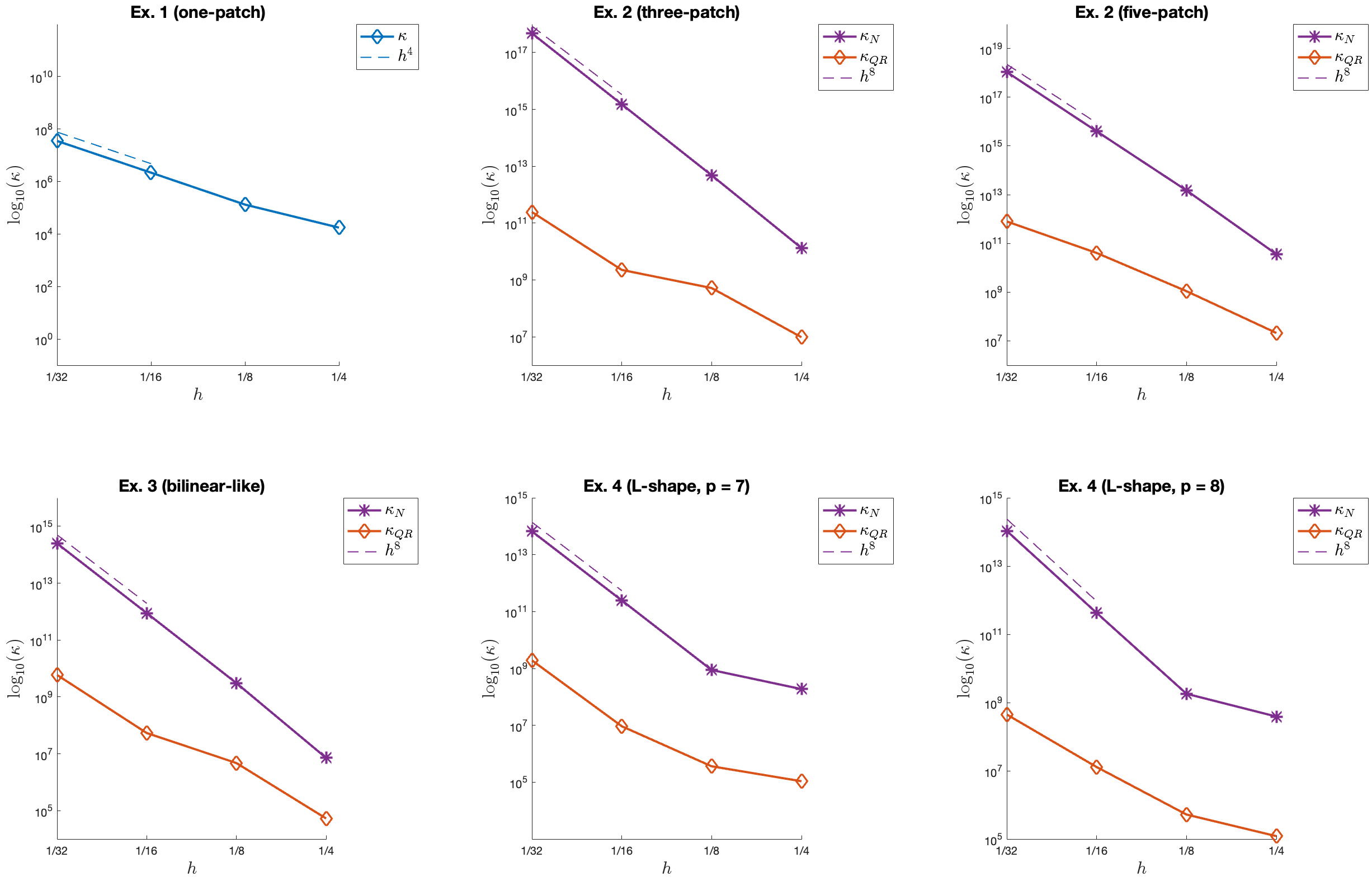}
    \caption{\new One patch example (top left): The condition numbers $\kappa$ for solving the square system. Multi-patch examples (top center -- bottom right): The condition numbers using the scaled normal system approach ($\kappa_N$) and QR decomposition ($\kappa_{QR}$) for solving the resulting overdetermined system.}
    \label{fig:conditionNumbers}
\end{figure}
\begin{table}[htb!]
\centering
{\new 
\begin{tabular}{|l|cccc|}
\hline
$\qquad$ Case $\,$ \textbackslash $\,$  $h$               & $1/4$            & $1/8$            & $1/16$             & $1/32$              \\ \hline
Ex.~1 (one-patch)
& $(625 , 625)$    & $(2025 , 2025)$  & $(7225 , 7225)$    & $(27225 , 27225)$   \\
Ex.~2 (three-patch)                 
& $(1795 , 1309)$  & $(5935 , 4969)$  & $(21415 , 19489)$  & $(81175 , 77329)$   \\
Ex.~2 (five-patch)                 
& $(2991 , 2166)$  & $(9891 , 8266)$  & $(35691 , 32466)$  & $(135291 , 128866)$ \\
Ex.~3 (bilinear-like) 
& $(2384 , 1720)$  & $(7904 , 6600)$  & $(28544 , 25960)$  & $(108224 , 103080)$ \\
Ex.~4 (L-shape, $p=7$) 
& $( 758 , 678)$   & $(2518 , 2358 )$ & $(9110 , 8790 )$   & $(34582 , 33942 )$  \\
Ex.~4 (L-shape, $p=8$) 
& $(1102 , 1014 )$ & $(3782 , 3606 )$ & $(13942 , 13590 )$ & $(53462 , 52758 )$  \\ \hline
\end{tabular}
\caption{Matrix dimensions of the collocation matrices obtained for the one-patch domain in Example~\ref{ex:onepatchdomain}, for the three-patch 
and five-patch domains in Example~\ref{ex:threepatchdomain}, for the bilinear-like four-patch domain in Example~\ref{ex:fourpatchdomainHole} and for the two-patch L-shape domain for $p=7$ and $p=8$ in Example~\ref{ex:twopatchdomainSmoothness3}.}
\label{tab:matrixSizes}
}
\end{table}
}

\section{Conclusion} \label{sec:Conclusion}

We presented a novel method for solving the biharmonic equation over planar bilinearly parameterized multi-patch domains in strong form, {\new and extended the presented technique on the basis of an example also to the case of bilinear-like $G^4$ multi-patch geometries, which allow the modeling of curved interfaces and boundaries}. Our developed isogeometric collocation approach is based on the use of the globally $C^4$-smooth isogeometric spline space~\cite{KaVi20b} 
with spline degree~$p=9$ to represent the $C^4$-smooth numerical solution of the biharmonic equation. For the choice of the collocation points, two different sets of points have been proposed, 
namely the Greville points, and a new family of superconvergent points. 
We tested both choices on a one-patch domain as well as on several instances of multi-patch domains. For both 
sets of collocation points, we numerically studied the convergence behavior with respect to the $L^2$-norm and with respect to (equivalents of) the $H^s$-seminorms, $1\leq s\leq 4$. For the considered odd spline degree~$p=9$, we observed in case of the Greville points a convergence of order~$\mathcal{O}(h^{p-3})$ for all (semi)norms, and in case of the superconvergent points an improved convergence of order~$\mathcal{O}(h^{p-2})$ for all (semi)norms except for the $H^4$-seminorm, where the order~$\mathcal{O}(h^{p-3})$ is optimal. 
 
We also studied on the basis of an example a first possibility to slightly reduce the required spline degree of the isogeometric multi-patch discretization space by reducing the regularity of the numerical solution as well. Another strategy and a first possible topic for future research could be to employ functions of degree $p=9$ just in the vicinity of edges and vertices, and to use functions with spline degree~$p=5$ in the interior of patches. 
A further open issue is the following. Since the number of collocation points in the multi-patch case is slightly larger than the dimension of the discretization space, the finding of a set of collocation points with the same cardinality as the dimension of the space would be of interest, too, in order to avoid the necessity of the least-squares method for solving the resulting linear system. 
In addition, other applications such as 
the Kirchhoff plate or Kirchhoff-Love shell problem, 
and the extension of our approach to multi-patch 
surfaces or to multi-patch volumes 
could be worth to study, too. 

\paragraph*{\bf Acknowledgment}

{\new The authors wish to thank the anonymous reviewers for their comments that helped to improve the paper.} M. Kapl has been partially supported by the Austrian Science Fund (FWF) through the project P~33023-N. V.~Vitrih has been partially supported by the 
Slovenian Research and Innovation Agency (research program P1-0404 and research projects N1-0296, J1-1715, N1-0210 and J1-4414). A.~Kosma\v c has been partially supported by the Slovenian Research and Innovation Agency (research program P1-0404, research project N1-0296 and Young Researchers Grant). This support is gratefully acknowledged.


\end{document}